\documentclass[12pt,dvips]{amsart}
\usepackage{amsfonts, amssymb, latexsym, epsfig}
\usepackage{euler, amsfonts, amssymb, latexsym, epsfig, epic}

\newtheorem{Theorem}{Theorem}[section] 

\newtheorem{Proposition}[Theorem]{Proposition} 
\newtheorem{Lemma}[Theorem]{Lemma}

\newtheorem{Corollary}[Theorem]{Corollary}
\newtheorem{Problem}[Theorem]{Problem}
\newtheorem{Definition-Lemma}[Theorem]{Definition-Lemma}
\theoremstyle{remark}
\newtheorem{Example}[Theorem]{Example}

\newcommand{\Aut}{\textrm{Aut}}

\theoremstyle{plain}

\newtheorem*{MainThm}{Main Theorem}

\newcommand{\cellsize}{13}
\newlength{\cellsz} 
\newsavebox{\cell}
\sbox{\cell}{\begin{picture}(\cellsize,\cellsize)
\put(0,0){\line(1,0){\cellsize}}
\put(0,0){\line(0,1){\cellsize}}
\put(\cellsize,0){\line(0,1){\cellsize}}
\put(0,\cellsize){\line(1,0){\cellsize}}
\end{picture}}
\newcommand\cellify[1]{\def\thearg{#1}\def\nothing{}%
\ifx\thearg\nothing
\vrule width0pt height\cellsz depth0pt\else
\hbox to 0pt{\usebox{\cell} \hss}\fi%
\vbox to \cellsz{
\vss
\hbox to \cellsz{\hss$#1$\hss}
\vss}}
\newcommand\tableau[1]{\vtop{\let\\\cr
\baselineskip -16000pt \lineskiplimit 16000pt \lineskip 0pt
\ialign{&\cellify{##}\cr#1\crcr}}}
\begin{document}
\pagestyle{plain}

\mbox{}\vspace{-2.0ex}
\title{A combinatorial rule for (co)minuscule Schubert calculus\vspace{-1ex}}

\author{Hugh Thomas}
\address{Department of Mathematics and Statistics, University of New Brunswick, Fredericton, New Brunswick, E3B 5A3, Canada }
\email{hugh@math.unb.ca}

\author{Alexander Yong}
\address{Department of Mathematics, University of Minnesota, Minneapolis, MN 55455, USA,
\indent{\itshape and} Department of Statistics/the Fields Institute, University of Toronto, Toronto, Ontario,  
        M5T 3J1, Canada}

\email{ayong@math.umn.edu, ayong@fields.utoronto.ca}

\date{September 7, 2006}

\maketitle
\vspace{-.2in}
\begin{abstract}\vspace{-2ex}
We prove a root system uniform, concise combinatorial rule
for Schubert calculus of \emph{minuscule} and \emph{cominuscule} flag manifolds~$G/P$ 
(the latter are also known as \emph{compact Hermitian symmetric spaces}).
We connect this geometry to the poset combinatorics of [Proctor '04], thereby
giving a generalization of the [Sch\"{u}tzenberger `77] \emph{jeu de taquin} 
formulation of the Littlewood-Richardson rule that computes the intersection
numbers of Grassmannian Schubert varieties. Our proof introduces
\emph{cominuscule recursions}, a general technique to relate the
numbers for different Lie types. A discussion about connections of our rule
to (geometric) representation theory is also briefly entertained.
\end{abstract}

\vspace{-3ex}
\tableofcontents

\section{Introduction}

\subsection{Overview} The goal of this paper is to introduce and prove 
a root-theoretically uniform 
generalization of the Littlewood-Richardson rule for 
intersection numbers of Schubert varieties in minuscule and
cominuscule flag varieties. 

Let $G$ denote a complex, connected, reductive (e.g., semisimple) 
Lie group. Fix a choice of Borel and opposite Borel 
subgroups~$B, B_{-}$ and maximal torus $T=B\cap B_{-}$. 
Let $W$ denote the Weyl group $N(T)/T$, 
$\Phi=\Phi^{+}\cup\Phi^{-}$ the 
ordering of the roots into positives and negatives, and 
$\Delta$ the base of simple roots. Choosing a parabolic
subgroup~$P$ canonically corresponds to a subset 
$\Delta_{P}\subseteq \Delta$;
let $W_P:=W_{\Delta_P}$ denote the associated parabolic subgroup of $W$. 
The {\bf generalized flag variety} $G/P$ is
a union of $B_{-}$-orbits whose closures 
${\mathcal X}_w:=\overline{B_{-}wP/P}$ with $wW_P\in W/W_P$
are the {\bf Schubert varieties}. The
Poincar\'{e} duals $\{\sigma_{w}\}$ of the Schubert varieties form the
{\bf Schubert basis} of the cohomology ring 
$H^{\star}(G/P)=H^{\star}(G/P; {\mathbb Q})$.

Among the simplest of the $G/P$'s are the projective spaces 
and Grassmannians. However, as our results help demonstrate,
the relative simplicity of their 
geometric and \linebreak representation-theoretic features is shared by the wider settings
of {\bf minuscule} and {\bf cominuscule flag varieties} 
(the latter are better known as {\bf compact Hermitian symmetric
spaces}). These are selected cases of~$G$ and its maximal parabolic
subgroup $P$; see the precise
definition in Section~2.1. (Actually there is 
little difference between the two settings.
We focus on the latter, explaining the adjustments for 
the former as necessary.)
These $G/P$'s and their
Schubert varieties are of significant and fundamental interest in geometry and
representation theory, see, e.g.,~\cite[Chapter~9]{Billey.Lakshmibai},~\cite{Kostant}
and the references therein, as well as, e.g.,~\cite{Perrin, Purbhoo.Sottile} 
for more recent work.

The generalities of $G/P$ specialize nicely to the
(co)minuscule cases. 
Associated to $G$ is the poset of positive roots 
$\Omega_{G}=(\Phi^{+},\prec)$ defined by the transitive closure of 
the covering relation $\alpha\prec \gamma$ if $\gamma-\alpha\in \Delta$.
For each (co)minuscule $G/P$ let $\Delta\setminus \Delta_{P}=\{\beta(P)\}$
be the simple root corresponding to $P$. We study the subposet 
\[\Lambda_{G/P}=\{\alpha\in \Phi^{+}: \mbox{$\alpha$ contains $\beta(P)$ in
its simple root expansion}\}\subseteq \Omega_{G}.\] 
The (co)minuscule hypothesis assures  that $\Lambda_{G/P}$ is self-dual and 
planar, see Section~2.2.

Moreover, rather than work with $W/W_P$-cosets directly, it is
possible in (co)minuscule cases to view the
Schubert basis as indexed by lower order 
ideals $\lambda\subseteq \Lambda_{G/P}$ 
(for a proof, see Proposition~\ref{prop:lower}).
We refer to these lower order ideals as {\bf (straight) shapes},
and we call their elements {\bf boxes}. 
Let ${\mathbb Y}_{G/P}$ be the lattice of 
these shapes, ordered by containment.

The {\bf Schubert intersection numbers} $\{c_{\lambda,\mu}^{\nu}(G/P)\}$ 
are defined by
\begin{equation}
\label{eqn:sin}
\sigma_{\lambda}\cdot \sigma_{\mu}=\sum_{\nu\in {\mathbb Y}_{G/P}}c_{\lambda,\mu}^{\nu}(G/P)\sigma_{\nu}.
\end{equation}
These numbers count points of intersection of generically translated
Schubert varieties 
and are therefore positive integers invariant under a natural 
$S_{3}$-action on the indices. 
 
It is a longstanding goal in combinatorial algebraic geometry to discover a
visibly positive combinatorial rule useful for understanding
the numbers $c_{\lambda,\mu}^{\nu}(G/P)$. 
Few cases have complete solutions or 
conjectures, even for $G=GL_{n}({\mathbb C})$.
The archetypal Grassmannian case 
is solved by the Littlewood-Richardson rule;
the first modern statement and proof is due to Sch\"{u}tzenberger~\cite{Schutzenberger} using the
combinatorics of jeu de taquin. 
See, e.g.,~\cite{buch:KLR, Coskun, knutson.tao:equivariant, Knutson.Yong, Kogan},
and the references therein, for variations on
generalized Littlewood-Richardson type formulas for 
$GL_{n}({\mathbb C})$-Schubert calculus.

It is a natural problem is to find such a rule 
for (co)minuscule flag varieties.

This paper extends the jeu de taquin formulation of
the Littlewood-Richardson rule to the Schubert calculus of
(co)minuscule flag varieties. It provides the first uniform 
generalization of the Littlewood-Richardson rule that involves both 
classical and exceptional Lie types (see, e.g.,~\cite{Purbhoo} for
earlier efforts in this direction). This suggests 
the potential to extend alternative frameworks for 
the Littlewood-Richardson rule and its consequences and/or generalizations
to the (co)minuscule setting, and beyond. 

In particular, our rule may be interpreted in terms of 
emerging and classical connections between Schubert calculus and 
(geometric) representation theory: for example, in connection to
the geometric Satake correspondence of Ginzburg, 
Mirkovi\'{c}-Vilonen and others, see, e.g.,~\cite{Mirkovic.Vilonen}, 
and separately, to Kostant's \cite{Kostant}
study of Lie algebra cohomology (see also~\cite{Belkale.Kumar}). See the 
remarks in Section~6.

\subsection{Statement of the main result}
If $\lambda\subseteq \nu$ are in ${\mathbb Y}_{G/P}$,
their set-theoretic difference is the {\bf skew shape} $\nu/\lambda$.
A {\bf standard filling} of a (skew) shape $\nu/\lambda$
is a bijective assignment
${\tt label}: \nu/\lambda\to \{1,2,\ldots,|\nu/\lambda|\}$
such that ${\tt label}(x)<{\tt label}(y)$ whenever $x\prec y$. The
result of this assignment is a {\bf standard tableau} $T$ of {\bf shape} 
$\nu/\lambda={\tt shape}(T)$. The set of all standard tableaux is
denoted ${\rm SYT}_{G/P}(\nu/\lambda)$. 
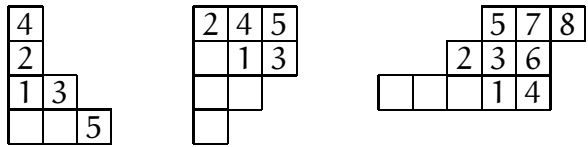
\begin{figure}[h]
\begin{picture}(170,40)
\put(0,30){$\tableau{{4}\\{2}\\{1}&{3}\\{\ }&{\ }&{5}}$}
\put(70,30){$\tableau{{2}&{4}&{5}\\{\ }&{1}&{3}\\{\ }&{\ }\\{\ }}$}
\put(140,30){$\tableau{&&&{5}&{7}&{8}\\&&{2}&{3}&{6}\\{\ }&{\ }&{\ }&{1}&{4}}$}
\end{picture}
\caption{Standard tableaux of shape $\nu/\lambda$ for types $A_{n-1}$, $C_n$ and $E_6$ respectively
(the empty boxes are those of $\lambda$); see Sections~2 and~3 for context.}
\end{figure}

Given $T\in {\rm SYT}_{G/P}(\nu/\lambda)$ we now present {\bf (co)minuscule
jeu de taquin}. Consider $x\in \Lambda_{G/P}$ that is not in $\nu/\lambda$,
maximal in $\prec$ subject to the condition that
it is below \emph{some} element of $\nu/\lambda$. 
We associate another standard tableau (of a different skew shape)
${\tt jdt}_{x}(T)$ arising from $T$ called the {\bf jeu de
taquin slide} of $T$ into $x$: Let $y$ be the box of $\nu/\lambda$
with the smallest label, among those that cover~$x$. Move ${\tt label}(y)$ to $x$, leaving
$y$ vacant. Look for boxes of $\nu/\lambda$ that cover $y$ and repeat the process,
moving into $y$ the smallest label available among those boxes covering it. 
The tableau ${\tt jdt}_{x}(T)$ is outputted when no such
moves are possible. (The result is a standard tableau; indeed, all the
intermediate tableaux are.) The {\bf rectification} of $T$
is the result of an iteration of {\bf jeu de taquin slides} until we 
terminate at a standard tableau ${\tt rectification}(T)$ of a (straight) shape.
\begin{figure}[h]
\begin{picture}(240,40)
\put(30,30){$\tableau{&&&&{7}\\&&&&{2}&{5}&{6}\\&&&{\ }&{1}&{4}\\{\ }&{\ }&{\ }&{\ }&{x }&{3}}$}
\put(180,30){$\tableau{&&&&{7}\\&&&&{5}&{6}&{\ }\\&&&{\ }&{2}&{4}\\{\ }&{\ }&{\ }&{\ }&{1}&{3}}$}
\put(0,20){$T=$}
\put(140,20){${\tt jdt}_{x}(T)=$}
\end{picture}
\caption{A standard tableau and a jeu de taquin slide, in type $E_7$}
\end{figure}
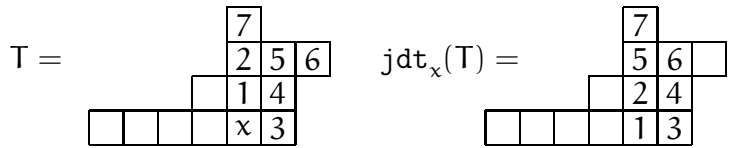

A novelty of this paper is the connection between (co)minuscule 
Schubert calculus and work of 
Proctor~\cite{Proctor}. That paper extends results of 
Sch\"{u}tzenberger~\cite{Schutzenberger},
Sagan~\cite{sagan} and Worley~\cite{worley}. It proves 
in the greater generality of 
``$d$-complete posets'' that the rectification 
is independent of the order of jeu de taquin slides; see Section~4.

Define the statistic ${\tt shortroots}$ on (skew) shapes to be the number of
boxes of $\nu/\lambda\subseteq \Lambda_{G/P}$ that are short roots. 
We are now ready to state our combinatorial rule. 
\begin{MainThm}
\label{thm:main}
Let $\lambda,\mu,\nu\in {\mathbb Y}_{G/P}$ and fix 
$T_{\mu}\in {\rm SYT}_{G/P}(\mu)$. In the minuscule case, 
the Schubert intersection number $c_{\lambda,\mu}^\nu(G/P)$ equals the number
of standard tableau of shape $\nu/\lambda$ whose 
rectification is~$T_{\mu}$; in the cominuscule case, 
multiply this by $2^{{\tt shortroots}(\nu/\lambda)-{\tt shortroots}(\mu)}$.
\end{MainThm}
The Main Theorem 
provides a root-theoretic generalization/reformulation of classical
theorems in the subject. Besides 
Sch\"utzenberger's rule~\cite{Schutzenberger} for Grassmannians, 
it generalizes the 
work of~\cite{Pragacz,worley} for isotropic Grassmannians. Moreover,
in the latter case, the power of $2$
that appears in the product rule for Schur $Q-$ polynomials 
is given a new interpretation, via the
${\tt shortroots}$ statistic. We emphasize that for the simply-laced 
root systems, the factor
$2^{{\tt shortroots}(\nu/\lambda)-{\tt shortroots}(\mu)}$ 
always equals~$2^{0-0}=1$.

In Section~2, we give preliminaries about 
(co)minuscule flag varieties and associated 
combinatorics. We give
examples of the Main Theorem in Section~3. Our proof, found in Sections~4--6,
is a collaboration of combinatorial and geometric ideas. 
There, we introduce the ideas of the ``infusion involution''
and ``cominuscule recursions''.
The latter are central to our (non-uniform) proof method of reducing
the difficult \emph{exceptional} Lie type cases to the \emph{classical} 
Lie type orthogonal Grassmannian case; this argument is based on the geometric
observation that certain Richardson varieties are isomorphic to Schubert
varieties in smaller cominuscule flag manifolds. 
In fact, these ideas are introduced in greater generality than needed in our proofs.
However, we believe they are interesting
in their own right, and we attempted to describe them in a natural
context.
We conclude in Section~7 with a collection of remarks and problems. 

The discovery of the Main Theorem exploited a number of computational
tools: John Stembridge's {\tt SF} and {\tt Coxeter/Weyl} packages for {\tt Maple}, 
Allen Knutson's algorithm~\cite{Knutson:alg} 
(as implemented in~\cite{Yong}).
In addition, we wrote the {\tt Maple} package 
{\tt Cominrule} to aid the reader in
exploring the properties of both the rule and (co)minuscule 
jeu de taquin.~\footnote{Available at the authors' websites.}

\section{(Co)minuscule flag varieties and their combinatorics}
\subsection{Definition and classification}
Our main source for background on (co)minuscule flag varieties 
is~\cite[Chapter~9]{Billey.Lakshmibai}.
For a
maximal parabolic subgroup $P$, interchangeably call it, its flag variety $G/P$ or
the root $\beta(P)\in\Delta$ (or, more properly, also 
the fundamental weight $\omega_{\beta(P)}$) 
{\bf cominuscule} if whenever $\beta(P)$ occurs in the simple root expansion
of $\gamma\in\Phi^{+}$, it does so with coefficient one. 
The cominuscule $G/P$'s have
been classified, see Table~\ref{table:classif}. In
each case, the cominuscule $\beta(P)\in\Delta$ are marked in the
Dynkin diagram. In case
of choice, selecting either one leads to a (possibly isomorphic) cominuscule
$G/P$. 

A maximal parabolic subgroup $P$, $G/P$ and $\beta(P)\in \Delta$
is {\bf minuscule} if the associated fundamental weight $\omega_{\beta(P)}$
satisfies $\langle\omega_{\beta(P)},\alpha^{\star}\rangle\leq 1$ for all 
$\alpha\in\Phi^{+}$ under the usual pairing between weights and coroots.
The classification of minuscule flag varieties \emph{almost} coincides 
with that of the cominuscules. In the conventions of Table~\ref{table:classif}, for the 
type $B_n$ minuscule case we select node $n$ rather than node $1$. 
This is the {\bf odd orthogonal Grassmannian} $OG(n, 2n+1)$, which is actually isomorphic
to the even orthogonal Grassmannian $OG(n+1, 2n+2)$. Consequently,
their Schubert intersection
numbers coincide. The type $C_n$ minuscule case corresponds to selecting node
$1$ rather than node $n$. 
\begin{center}
\begin{table}
\begin{tabular}{||c|c|c||}\hline
\emph{Root system} & \emph{Dynkin Diagram} & \emph{Nomenclature}\\\hline \hline
$A_n$ & 
\setlength{\unitlength}{3mm} 
\begin{picture}(11,3)
\multiput(0,1.5)(2,0){6}{$\circ$}
\multiput(0.55,1.85)(2,0){5}{\line(1,0){1.55}}
\put(6,1.5){$\bullet$}
\put(0,0){$1$}
\put(2,0){$2$}
\put(3.5,0){$\cdots$}
\put(6,0){$k$}
\put(7.5,0){$\cdots$}
\put(10,0){$n$}
\end{picture}
& Grassmannian $Gr(k,{\mathbb C}^{n})$\\ \hline
$B_n$ & 
\setlength{\unitlength}{3mm} \begin{picture}(11,3)
\multiput(0,1.5)(2,0){6}{$\circ$}
\multiput(0.55,1.85)(2,0){4}{\line(1,0){1.55}}
\multiput(8.55,1.75)(0,.2){2}{\line(1,0){1.55}} \put(8.75,1.52){$>$}
\put(0,1.5){$\bullet$}
\put(0,0){$1$}
\put(2,0){$2$}
\put(4,0){$\cdots$}
\put(7,0){$\cdots$}
\put(10,0){$n$}\end{picture}
& Odd dimensional quadric ${\mathbb Q}^{2n-1}$\\ \hline
$C_n, n\geq 3$ & 
\setlength{\unitlength}{3mm} \begin{picture}(11,3)
\multiput(0,1.5)(2,0){6}{$\circ$}
\multiput(0.55,1.85)(2,0){4}{\line(1,0){1.55}}
\multiput(8.55,1.75)(0,.2){2}{\line(1,0){1.55}} \put(8.85,1.53){$<$}
\put(10,1.5){$\bullet$}
\put(0,0){$1$}
\put(2,0){$2$}
\put(4,0){$\cdots$}
\put(7,0){$\cdots$}
\put(10,0){$n$}
\end{picture}
& Lagrangian Grassmannian $LG(n,2n)$\\ \hline
$D_n, n\geq 4$ & 
$\begin{array}{c}
\setlength{\unitlength}{2.9mm} \begin{picture}(11,3.5)
\multiput(0,1.6)(2,0){5}{$\circ$}
\multiput(0.55,2)(2,0){4}{\line(1,0){1.55}} 
\put(8.5,1.95){\line(2,-1){1.55}}
\put(8.5,1.95){\line(2,1){1.55}}
\put(10,2.5){$\circ$}
\put(10,0.7){$\circ$}
\put(0,1.6){$\bullet$}
\put(0,0){$1$}
\put(2,0){$2$}
\put(4,0){$\cdots$}
\put(7,0){$\cdots$}
\put(9.1,0){$n\!-\!1$}
\put(11, 2.3){$n$}\end{picture}
\\
\setlength{\unitlength}{2.9mm} \begin{picture}(11,3.5)
\multiput(0,1.6)(2,0){5}{$\circ$}
\multiput(0.55,2)(2,0){4}{\line(1,0){1.55}} 
\put(8.5,1.95){\line(2,-1){1.55}}
\put(8.5,1.95){\line(2,1){1.55}}
\put(10,2.5){$\circ$}
\put(10,0.7){$\circ$}
\put(0,0){$1$}
\put(2,0){$2$}
\put(4,0){$\cdots$}
\put(7,0){$\cdots$}
\put(9.1,0){$n\!-\!1$}
\put(11, 2.3){$n$}

\put(10,2.45){$\bullet$}
\put(10,0.75){$\bullet$}

\end{picture}
\end{array}$

& $\begin{array}{c}
\mbox{Even dimensional quadric ${\mathbb Q}^{2n-2}$}\\
\\
\mbox{Orthogonal Grassmannian $OG(n+1,2n+2)$}\\
(\mbox{for either choice of one of the nodes $n-1$ or $n$})
\end{array}$
\\\hline         
$E_6$ & 
\setlength{\unitlength}{3mm}
\begin{picture}(9,3.6)
\multiput(0,0.5)(2,0){5}{$\circ$}
\multiput(0.55,0.95)(2,0){4}{\line(1,0){1.6}} 
\put(0,0.5){$\bullet$}
\put(8,0.5){$\bullet$}
\put(4,2.6){$\circ$}
\put(4.35,1.2){\line(0,1){1.5}}
\put(0,-.6){$1$}
\put(2,-0.6){$3$}
\put(4,-.6){$4$}
\put(6,-.6){$5$}
\put(5,2.5){$2$}
\put(8,-.6){$6$}
\end{picture}
& 
$\begin{array}{c}
\mbox{Cayley Plane ${\mathbb O}{\mathbb P}^2$}\\ 
\mbox{(for either choice of one of the nodes $1$ or $6$)}\\ 
\end{array}$\\ \hline
$E_7$ & 
\setlength{\unitlength}{3mm}
\begin{picture}(11,4)
\put(0,0.9){$\bullet$}
\multiput(2,0.9)(2,0){5}{$\circ$}
\multiput(0.55,1.35)(2,0){5}{\line(1,0){1.6}} 
\put(10,0.9){$\circ$}
\put(6,3){$\circ$}
\put(6.35,1.6){\line(0,1){1.5}}
\put(0,-.2){$1$}
\put(2,-0.2){$3$}
\put(4,-.2){$4$}
\put(6,-.2){$5$}
\put(7,2.9){$2$}
\put(8,-.2){$6$}
\put(10,-.2){$7$}
\end{picture}
& 
$\mbox{(Unnamed) } G_{\omega}({\mathbb O}^3, {\mathbb O}^6)$
\\ \hline
\end{tabular}
\caption{\label{table:classif} Classification of cominuscule $G/P$'s}
\end{table}
\end{center}
\vspace{-.4in}
\subsection{More specifics about cominuscule $\Lambda_{G/P}$}
We
give two pictures of $\Lambda_{G/P}$: ``thin'' and
``fattened'' versions where elements of $\Lambda_{G/P}$ (i.e., boxes) are
represented by either dots ``$\bullet$'' or squares ``$\square$''.
The former are more convenient for illustrating the poset relations
(smaller elements are lower in the diagram), and
$\Lambda_{G/P}$'s position inside $\Omega_{G}$.
The latter are more convenient for manipulations in Section~3 and~4; for this reason, the 
fattened $\Lambda_{G/P}$ is rotated 45 degrees clockwise relative to the thin $\Lambda_{G/P}$.

The self-duality of $\Lambda_{G/P}$ is as follows. Let $u_0$ be the maximal length element of~$W_P$.
Note that $u_0$ preserves the subset 
$\Lambda_{G/P}$: the positive roots it makes negative
are exactly those not in $\Lambda_{G/P}$, and thus
if $u_0$ moved a root in $\Lambda_{G/P}$ outside, applying $u_0$ twice
would send that root negative, contradicting the fact $u_0^2=1$.
Let ${\tt rotate}$ denote this involution on $\Lambda_{G/P}$. 
In particular, this sends shapes to \emph{upper} order ideals of $\Lambda_{G/P}$ and conversely.

We summarize features of each $\Lambda_{G/P}$ in Table~\ref{table:Lambda}.
We proceed to analyze the specifics in each of the Lie types. 
Inside the fattened depictions, we draw a sample shape 
$\lambda=(\lambda_1,\lambda_2,\ldots\ )$ where $\lambda_i$
is the number of boxes in column~$i$, as read from left to right.

\begin{table}
\begin{tabular}{||c|c|c|c||}\hline
$G/P$ & \emph{$\Lambda_{G/P}$ description} & \emph{$\#$ boxes } & \emph{Short root boxes}\\\hline\hline
$Gr(k,{\mathbb C^n})$ & $k\times(n-k)$ rectangle & $k(n-k)$ & none \\ \hline
${\mathbb Q}^{2n-1}$ & $(2n-1)$-row & $2n-1$ & middle box \\ \hline
$LG(n,2n)$ & $n$-step staircase & ${n+1 \choose 2}$ & all \emph{non}-anti-diagonal boxes \\ \hline
${\mathbb Q}^{2n-2}$ & double tailed diamond & $2n-2$ & none\\ \hline
$OG(n+1,2n+2)$ & $n-1$-step staircase & ${n\choose 2}$ & none \\ \hline
${\mathbb O}{\mathbb P}^{2}$ & irregular & 16 & none\\ \hline
$G_{\omega}({\mathbb O}^3, {\mathbb O}^6)$ & irregular & 27 & none \\ \hline
\end{tabular}
\caption{\label{table:Lambda} Summary of facts about cominuscule $\Lambda_{G/P}$}
\end{table}

\noindent\emph{Type $A_{n-1}$:} 
\begin{figure}[h]
\begin{center}
\epsfig{file=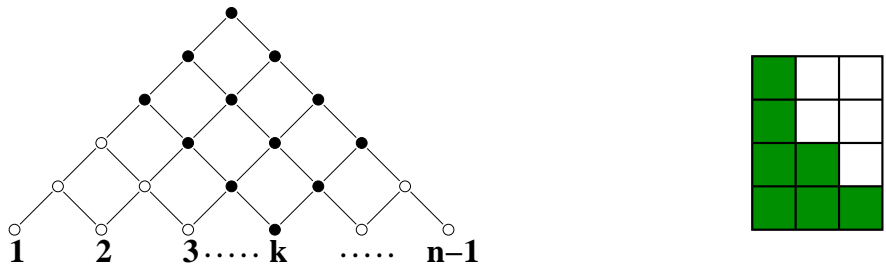, height=3.4cm}
\end{center}
\caption{\label{fig:An} $\Lambda_{Gr(k,{\mathbb C}^n)}$, $\Omega_{GL_{n}}({\mathbb C})$ and the shape $\nu=(4,2,1)$}
\end{figure}
Figure~\ref{fig:An} depicts the case $k=4$ and $n=7$. 
Lower order ideals correspond to 
Young shapes (partitions) drawn in ``conjugate French notation''.
The ${\tt rotate}$ involution is ``rotate by $180$ degrees''.
So, for example, ${\tt rotate}(\nu)$ is the complement $(3,2)^c$
in $\Lambda_{Gr(k,{\mathbb C}^n)}$.

\begin{figure}[h]
\begin{center}
\epsfig{file=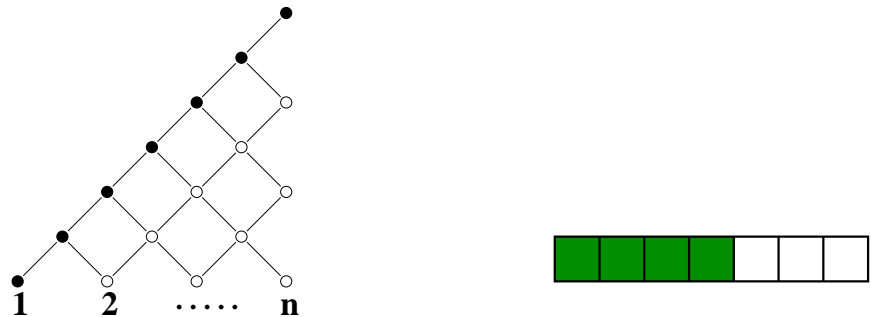, height=4.1cm}
\end{center}
\caption{\label{fig:Bn} $\Lambda_{{\mathbb Q}^{2n-1}}$, $\Omega_{SO_{2n+1}({\mathbb C})}$ and
the shape $\nu=(1^4)=(1,1,1,1)$}
\end{figure}

\noindent\emph{Type $B_n$:} Classes are indexed by a single row of some 
length~$j$,
denoted $(1^j)=(1,1,\ldots,1)$. 
Again ${\tt rotate}$ is ``rotate by $180$ degrees''.
See Figure~\ref{fig:Bn}. Here, ${\tt rotate}(1^4)=(1^3)^c$.

\begin{figure}[h]
\begin{center}
\epsfig{file=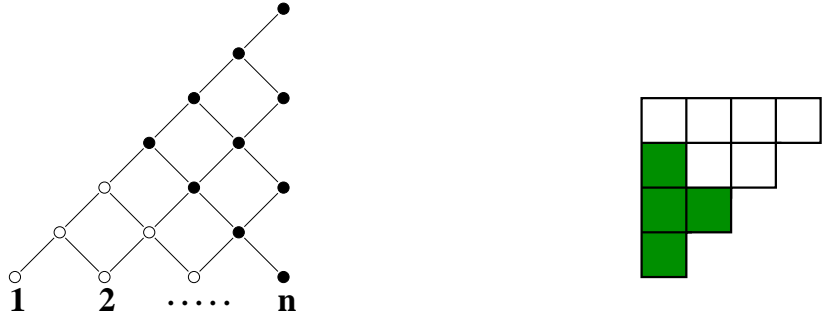, height=4.1cm}
\end{center}
\caption{\label{fig:Cn} $\Lambda_{LG(n,2n)}$, $\Omega_{Sp_{2n}({\mathbb C})}$ and
the shape $\nu=(3,1)$}
\end{figure}

\noindent\emph{Type $C_n$:} Shapes are {\bf strict partitions} 
$\lambda=(\lambda_1>\lambda_2>\ldots \ )$ contained inside the 
``staircase''. The boxes \emph{not} on the {\bf anti-diagonal} of the
fattened diagram (i.e., those boxes not 
directly above ``$n$'' in the thin diagram) 
correspond to short roots. Here {\tt rotate} corresponds to flipping
across the diagonal line of symmetry in the fattened depiction.
See Figure~\ref{fig:Cn}, where, e.g., the shape
$(3,1)$ involves two short roots, the shape $(4,1)$ would involve three
short roots. Here ${\tt rotate}(3,1)=(4,2)^c$.

\noindent\emph{Type $D_n$:} 
The Orthogonal Grassmannian $OG(n+1,2n+2)$ comes from either choosing
node $n-1$ or $n$. In either case, $\Lambda_{OG(n+1,2n+2)}$ 
is isomorphic (as a poset) to $\Lambda_{LG(n-1,2n-2)}$ (type $C_{n-1}$) drawn above. 
The shapes are shifted shapes, as in type $C_{n-1}$ above.
However, this time, none of the boxes correspond to short roots. 

\begin{figure}[h]
\begin{center}
\epsfig{file=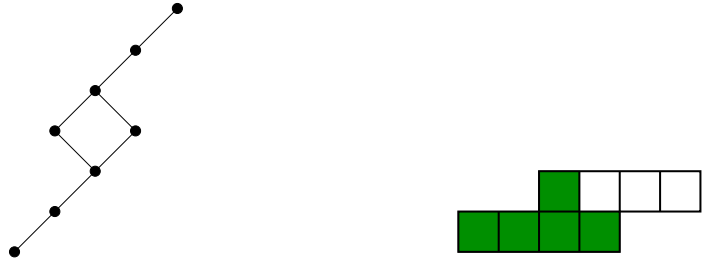, height=3.4cm}
\end{center}
\caption{\label{fig:Dnquadric} $\Lambda_{{\mathbb Q}^{2n-2}}$ and
the shape $\nu=(1,1,2,1)$}
\end{figure}

For the case of the even dimensional quadrics, see 
Figure~\ref{fig:Dnquadric}. (For the Dynkin diagrams with ``forks'',
the $\Omega_{G}$ becomes complicated to draw, so we omit them.) In this case
the visualization of {\tt rotate} depends on the parity of $n$. For the
$n$ is odd case it is the $180$ degree rotation, whereas 
if $n$ is even the same is true
except that the middle nodes stay fixed. For instance, if $n=5$, 
${\tt rotate}(1^4)$ is $(1^4)^c$ whereas when $n=6$,
${\tt rotate}(1^5)$ is $(1,1,1,2)^c$.

\begin{figure}[h]
\begin{center}
\epsfig{file=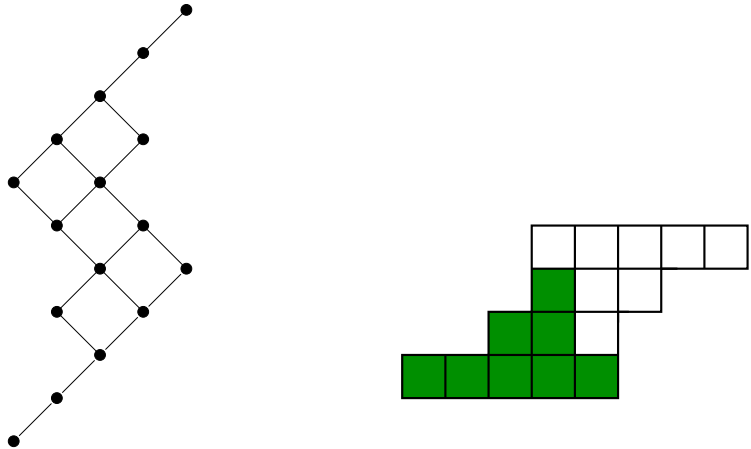, height=5.9cm}
\end{center}
\caption{\label{fig:E6} $\Lambda_{{\mathbb O}{\mathbb P}^2}$ and
the shape $\nu=(1,1,2,3,1)$}
\end{figure}

\begin{figure}[h]
\begin{center}
\epsfig{file=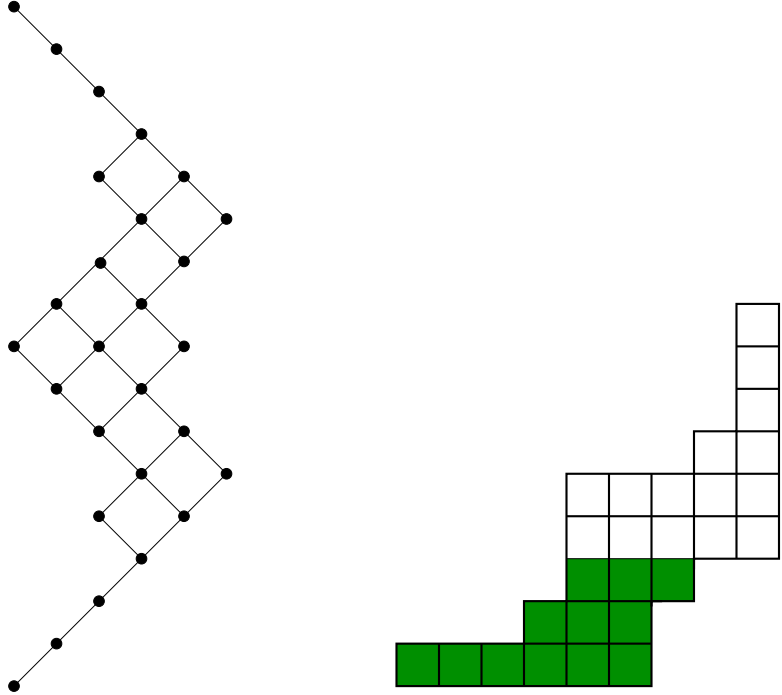, height=9.2cm}
\end{center}
\caption{\label{fig:E7} $\Lambda_{G_{\omega}({\mathbb O}^3, {\mathbb O}^6)}$ and
the shape $\nu=(1,1,1,2,3,3,1)$}
\end{figure}

\noindent
\emph{Types $E_6$ and $E_7$:} ${\tt rotate}$ is
again 180 degree rotation. See Figures~\ref{fig:E6} and~\ref{fig:E7}.
In the depicted $E_6$ case, ${\tt rotate}(\nu)=\nu^c$, whereas in the $E_7$ case
${\tt rotate}(\nu)=(1,1,1,2,5,5)^c$.

\subsection{Minuscule $\Lambda_{G/P}$}

As mentioned, the minuscule cases coincide with the cominuscules
except in types~$B$ and~$C$, which we discuss now.

\noindent
\emph{Type $B_n$:} This minuscule flag variety is isomorphic to
$OG(n+1,2n+2)$; its ambient poset is the same $\Lambda_{OG(n+1,2n+2)}$,
see Figure~\ref{fig:Cn}.
(Which nodes in the poset correspond to short roots is different, but,
because we are now in the minuscule case of the Main Theorem, this
is irrelevant.) Note that $u_0$ in this case now acts differently on
$\Lambda_{OG(n, 2n+1)}$, so in this case we declare that ${\tt rotate}$
acts as it does on ${\Lambda}_{OG(n+1,2n+2)}\cong 
{\Lambda}_{LG(n-1,2n-2)}$.

\noindent
\emph{Type $C_n$:} This is the projective space ${\mathbb P}^{2n-1}$. Here
$\Lambda_{G/P}\cong\Lambda_{{\mathbb Q}^{2n-1}}$ 
with the same shapes,
see Figure~\ref{fig:Bn}.
(Again, which nodes correspond to short and long roots differ, which
is irrelevant to the Main Theorem.) 
 
\subsection{Shapes index the Schubert basis\label{section:shapesindex}}
For $w\in W$ the {\bf inversion set} is
${\mathcal I}(w)=\{\alpha\in \Phi^{+}: w\cdot\alpha\in \Phi^{-}\}$,
using the standard action of $W$ on $\Phi$. 
Consider an irreducible rank two root system $\Phi_{\{\eta,\gamma\}}\subseteq \Phi$.  
It inherits from $\Phi$ a decomposition into positives and negatives, 
with $\eta,\gamma\in\Phi^{+}$ as its simple roots. 
If $\eta$ and $\gamma$ are of different lengths, let
$\gamma$ be the shorter of them. 
Order the positive roots of $\Phi_{\{\eta,\gamma\}}$ 
\begin{equation}
\label{eqn:posroots2}
(\eta, \eta+\gamma,\gamma) \mbox{ or } (\eta,\eta+\gamma,\eta+2\gamma,\gamma)
\end{equation}
depending on whether the root system is $A_2$ or $B_2$, respectively  
($G_2$ has no cominuscule simple roots). 
A subset
$S\subseteq \Phi^+$ is called {\bf biconvex} if $S\cap \Phi_{\{\eta,\gamma\}}$ 
is either a
beginning or ending subset of the positive roots, for all 
$\Phi_{\eta,\gamma}\subseteq \Phi^{+}$,
with respect to (\ref{eqn:posroots2}).   
It is
known that $S={\mathcal I}(w)$ for some $w\in W$ if and only if $S$ is biconvex \cite{Bj,Ed}. 

The natural projection $G/B\twoheadrightarrow G/P$ 
induces an inclusion 
$H^{\star}(G/P)\hookrightarrow H^{\star}(G/B)$ 
sending $\sigma_{wW_P}\in H^{\star}(G/P)$
to $\sigma_{w^P}\in H^{\star}(G/B)$, where $w^P\in W$
is the minimal length representative of $wW_P$. Alternatively, 
$\sigma_{w}\in H^{\star}(G/B)$
appears as the image of a Schubert class under the
projection if and only if the {\bf descents} 
of $w$ are a subset of $\Delta\setminus\Delta_{P}$, i.e., $\ell(ws_{\beta})<\ell(w)$ 
only when $\beta\in \Delta\setminus\Delta_{P}$. Here $\ell(w)$ denotes
the {\bf Coxeter length}, the minimal length of an 
expression for $w$ in terms of the 
{\bf simple reflections} $s_\beta$.
When $\Delta\setminus\Delta_{P}=\{\beta(P)\}$ we call such a $w$ 
{\bf Grassmannian at} $\beta(P)$. Equivalently, these are the elements of $W^P$, the
minimal length coset representatives of $W/W_P$.

The remaining facts in this section are essentially well-known. 
We include a proof of the following since it is basic
to this paper.
\begin{Proposition}
\label{prop:lower}
Let $G/P$ be a cominuscule flag variety. The Schubert 
classes in $H^{\star}(G/P)$
are in bijection with $\lambda\in {\mathbb Y}_{G/P}$. 
Specifically, for each
$w\in W$ that is Grassmannian at~$\beta(P)$, the inversion set 
${\mathcal I}(w)$ is a lower
order ideal in $\Lambda_{G/P}$, and conversely every lower order ideal in
$\Lambda_{G/P}$ is the inversion set of a $w\in W$ which is Grassmannian
at $\beta(P)$.  
\end{Proposition}
\begin{proof}
Let $\lambda\in {\mathbb Y}_{G/P}$. 
Let $\Phi_{\{\eta,\gamma\}}$ be as above. If 
$\lambda\cap \Phi_{\{\eta,\gamma\}}\neq \emptyset$,
then by (\ref{eqn:posroots2}) either $\eta$ or 
$\gamma$ involve $\beta(P)$
in their simple root expansions.

By the cominuscule
assumption, any $\alpha\in \Phi^{+}$ involving $\beta(P)$ 
does so with coefficient one. From this and
(\ref{eqn:posroots2}) combined we deduce
that $\beta(P)$ appears in only one of $\eta$ or $\gamma$, and 
in the nonsimply-laced case, that $\beta(P)$ appears in $\eta$.
By symmetry, we may assume that $\beta(P)$ also appears in $\eta$ in the
simply-laced case. Thus $\lambda\cap {\Phi}_{\{\eta,\gamma\}}$ 
contains $\eta$ and since $\lambda$ is a lower order ideal, it
is a beginning subset of the positive roots. Hence $\lambda$ is
biconvex and $\lambda={\mathcal I}(w)$ for some $w\in W$. But the descents
of $w$ are $\Delta\cap {\mathcal I}(w)=\{\beta(P)\}$. So $w$ is
Grassmannian at $\beta(P)$ as desired.

Conversely, let $w$ be Grassmannian at $\beta(P)$. 
Suppose $\alpha\in {\mathcal I}(w)$ but $\beta(P)$ 
does not occur in $\alpha$. If $\alpha\not\in\Delta$ then write 
$\alpha=\gamma_1+\gamma_2$
where $\gamma_1,\gamma_2\in\Phi^{+}$. Notice at least one of $\gamma_1$ or
$\gamma_2$ is in ${\mathcal I}(w)$ (otherwise $\alpha$ would not be, a contradiction). 
Thus inductively, we reduce to the case
$\alpha\in\Delta$ anyway. Thus $\alpha$ is a descent of $w$, contradicting
our assumption that $w$ is Grassmannian at $\beta(P)$.   
So $\alpha$ involves $\beta(P)$, and
since $\beta(P)$ is cominuscule, it does so with coefficient one. Hence 
${\mathcal I}(w)\subseteq \Lambda_{G/P}$.

Now, suppose $\gamma\in {\mathcal I}(w)$ and $\delta\prec \gamma$ in
$\Lambda_{G/P}$. We may assume
$\gamma-\delta=\rho\in\Delta$. Note $\rho\neq\beta(P)$.
Since $\rho\not\in {\mathcal I}(w)$
but $\gamma\in {\mathcal I}(w)$, then $w\cdot\delta=w\cdot\gamma-w\cdot\rho\in\Phi^{-}$.
Thus $\delta\in {\mathcal I}(w)$, and so ${\mathcal I}(w)$ is a lower order ideal. 
\end{proof}

For brevity, we omit the proof of the next fact.
\begin{Lemma}\label{lemma:bruhatweak}
Let $G/P$ be a cominuscule flag variety and 
$\lambda,\nu\in {\mathbb Y}_{G/P}$. 
Then $\lambda\subseteq \nu$ if and only if 
$uW_P$ is smaller than $vW_P$ in the Bruhat order $W/W_P$
where $\lambda={\mathcal I}(u)$ and $\nu={\mathcal I}(v)$,
under the correspondence of Proposition~\ref{prop:lower}.
\end{Lemma}

\begin{Corollary}
\label{cor:overlap}
\begin{itemize}
\item[(a)] If $\lambda\not\subseteq\nu$ then 
$c_{\lambda,\mu}^{\nu}(G/P)=0$ for all shapes $\mu\subseteq \Lambda_{G/P}$.
\item[(b)] If $|\lambda|+|\mu|=|\Lambda_{G/P}|$
then $\sigma_{\lambda}\cdot\sigma_{\mu}=\left\{
\begin{array}{cc}
\sigma_{\Lambda_{G/P}} & \mbox{if $\ \lambda={\tt rotate}(\mu^c)$,}\\
0 & \mbox{otherwise,}
\end{array}\right.$
where $\mu^c$ is the complement of $\mu\subseteq \Lambda_{G/P}$.
\item[(c)] If $\lambda\cap {\tt rotate}(\mu)\neq \emptyset$ then
$\sigma_{\lambda}\cdot\sigma_{\mu}=0$.
\end{itemize}
\end{Corollary}
\begin{proof}
By Lemma~\ref{lemma:bruhatweak} and Proposition~\ref{prop:lower}, 
combined with the discussion of the injection 
$H^{\star}(G/P)\hookrightarrow H^{\star}(G/B)$
given before Proposition~\ref{prop:lower}, the assertions become well-known
facts about the Schubert intersection numbers on $G/B$.
\end{proof}

\section{Examples of the combinatorial rule \label{sect:exa}}

In the examples, $T_{\mu}$ is the 
``consecutive'' standard tableau, i.e.,  
the one with $1,2,3,\ldots,\mu_1$ labeling the first column, followed
by $\mu_1 +1, \mu_1 +2,\ldots$ labeling the second column etc. 

\noindent\emph{The Grassmannian:}
Let us do the computation
$c_{(3,1),(2,1)}^{(4,2,1)}(Gr(4, {\mathbb C}^7))=2$ (type $A_7$).
Here $\nu/\lambda = (4,2,1)/(3,1)$. Of the 6 standard tableaux,
only two rectify to $T_{2,1}$:
\[\tableau{{2}\\{\ }\\{\ }&{1}\\{\ }&{\ }&{3}}\ , \ \  \tableau{{2}\\{\ }\\{\ }&{3}\\{\ }&{\ }&{1}}
\mapsto \ \
\tableau{{\ }\\{\ }\\{2 }&{\ }\\{1 }&{3 }&{\ }}.\]

\noindent\emph{The isotropic Grassmannians:} First we compute 
$c_{(2,1),(2,1)}^{(4,2)}(LG(4,8))$ (type~$C_4$).
Here\linebreak $\nu/\lambda = (4,2)/(2,1)$ and ${\rm shortroots}(\nu/\lambda)=3$,
while ${\rm shortroots}(\mu)=1$. There are two standard tableau of shape
$\nu/\lambda$, but only one rectifies to $T_{2,1}$:
\[\tableau{{2}\\{1}&{3}\\{\ }&{\ }\\{\ }} \ \mapsto \ \ 
\tableau{{\ }\\{\ }&{\ }\\{2 }&{3 }\\{1 }}\]
Hence $c_{(2,1),(2,1)}^{(3,1)}(LG(4,8))=1\cdot 2^{3-1}=4$.
The same analysis shows $c_{(2,1),(2,1)}^{(4,2)}(OG(6,12))=1$ (type $D_{5}$),
since $\Lambda_{OG(6,12)}\cong \Lambda_{LG(4,8)}$
but now there are no short roots.

\noindent\emph{The quadrics:}
First let us compute $c_{(1^2),(1^2)}^{(1^4)}({\mathbb Q}^7)$ (type $B_4$). 
Here $\nu/\lambda=(1^4)/(1^2)$ and ${\rm shortroots}(\nu/\lambda)=1$
since this skew shape involves the middle box of $\Lambda_{{\mathbb Q}^{7}}$,
while \linebreak
${\rm shortroots}(\mu)=0$. We have 
$c_{(1^2),(1^2)}^{(1^4)}({\mathbb Q}^7)=1\cdot 2^{1-0}=2$ since:
\[\tableau{{\ }&{\ }&{1}&{2}} \ \mapsto \ \tableau{{1}&{2}&{\ }&{\ }}\]

The even dimensional 
quadrics have a quirky dependency on the parity of $n$.
When $n=5$, we have $c_{(1,1,2),(1,1,2)}^{\Lambda_{{\mathbb Q}^8}}(\Lambda_{{\mathbb Q}^8})=1$ as witnessed
by
\[\tableau{&&{\ }&{2}&{3}&{4}\\{\ }&{\ }&{\ }&{1}}\mapsto
\tableau{&&{4}&{\ }&{\ }&{\ }\\{1 }&{2 }&{3 }&{\ }},\]
whereas $c_{(1,1,2),(1^4)}^{\Lambda_{{\mathbb Q}^8}}(\Lambda_{{\mathbb Q}^8})=0$.
The similar $n=6$ computation gives
$c_{(1,1,1,2), (1^5)}^{\Lambda_{{\mathbb Q}^{10}}({\mathbb Q}^{10}})=1$ because
\[\tableau{&&&{\ }&{2}&{3}&{4}&{5}\\{\ }&{\ }&{\ }&{\ }&{1}}\mapsto
\tableau{&&&{\ }&{\ }&{\ }&{\ }&{\ }\\{1 }&{2 }&{3 }&{4 }&{5}}.\]
Thus cominuscule jeu de taquin properly detects the subtle  definition
of ${\tt rotate}$ in these cases: these calculations agree with 
Corollary~\ref{cor:overlap}(b).

\noindent\emph{The Cayley plane:}
We compute 
$c_{(1,1,2,1,1), (1,1,2,2,1)}^{(1,1,2,4,4,1)}({\mathbb O}{\mathbb P}^2)=2$
as shown by 
\[\tableau{&&&{4}&{6}\\&&&{2}&{5}&{7}\\&&{\ }&{1}&{3}\\{\ }&{\ }&{\ } & {\ }& {\ }}\ , \ 
\tableau{&&&{6}&{7}\\&&&{2}&{4}&{5}\\&&{\ }&{1}&{3}\\{\ }&{\ }&{\ } & {\ }& {\ }}\ \mapsto \ 
\tableau{&&&{\ }&{\ }\\&&&{\ }&{\ }&{\ }\\&&{4 }&{6}&{\ }\\{1 }&{2 }&{3 } & {5 }& {7 }}.\]

\noindent\emph{$G_{\omega}({\mathbb O}^3, {\mathbb O}^6)$:}
Finally, 
$c_{(1,1,1,2,5,3), (1,1,1,2,1)}^{(1,1,1,2,5,5,2,1,1)}(G_{\omega}({\mathbb O}^3, {\mathbb O}^6))=4$
since
\[\tableau{&&&&{\ }&{3}\\&&&&{\ }&{2}&{4}&{5}&{6}\\&&&&{\ }&{\ }&{1}\\
&&&{\ }&{\ }&{\ }\\{\ }&{\ }&{\ }&{\ }&{\ }&{\ }} \ , \!\!\!\!\!\!\!\!\!\!\!\!\!\!\!\!\!\!\!\!
\tableau{&&&&{\ }&{3}\\&&&&{\ }&{1}&{4}&{5}&{6}\\&&&&{\ }&{\ }&{2}\\
&&&{\ }&{\ }&{\ }\\{\ }&{\ }&{\ }&{\ }&{\ }&{\ }}\ , \!\!\!\!\!\!\!\!\!\!\!\!\!\!\!\!\!\!\!\!
\tableau{&&&&{\ }&{5}\\&&&&{\ }&{1}&{3}&{4}&{6}\\&&&&{\ }&{\ }&{2}\\
&&&{\ }&{\ }&{\ }\\{\ }&{\ }&{\ }&{\ }&{\ }&{\ }} \ , \!\!\!\!\!\!\!\!\!\!\!\!\!\!\!\!\!\!\!\!
\tableau{&&&&{\ }&{6}\\&&&&{\ }&{1}&{3}&{4}&{5}\\&&&&{\ }&{\ }&{2}\\
&&&{\ }&{\ }&{\ }\\{\ }&{\ }&{\ }&{\ }&{\ }&{\ }}
\mapsto \!\!\!\!\!\!\!\!\!\!\!\!\!\!\!\!\!\!\!\!\!\!
\tableau{&&&&{\ }&{\ }\\&&&&{\ }&{\ }&{\ }&{\ }&{\ }\\&&&&{\ }&{\ }&{\ }\\
&&&{5 }&{\ }&{\ }\\{1 }&{2 }&{3 }&{4 }&{6 }&{\ }}.
\] 

\section{Jeu de taquin methods}

The basic result used is due to Robert Proctor~\cite{Proctor}. We 
develop consequences of it for our purposes.
\begin{Theorem}\cite{Proctor}
\label{thm:proctor}
Let $\nu/\lambda\subseteq \Lambda_{G/P}$ be a skew shape and 
$T\in {\rm SYT}_{G/P}(\nu/\lambda)$. Then the procedure 
${\tt rectification}(T)$ is
independent of the order of applying jeu de taquin slides.
\end{Theorem}
In~\cite{Proctor}, this theorem is proved for the 
more general class of ``$d$-complete posets''. 
Mark Haiman has raised the question of assigning a geometric
context to these posets.

\subsection{Reversing}
Jeu de taquin is reversible. Given 
$T\in {\rm SYT}_{G/P}(\nu/\lambda)$, consider $x\in\Lambda_{G/P}$ not in
$\nu/\lambda$ but \emph{minimal} in $\prec$ subject to being
above some element of $\nu/\lambda$. The {\bf reverse jeu de taquin
slide} ${\tt revjdt}_{x}(T)$ of $T$ into $x$ is defined similarly
to a jeu de taquin slide, except we move into $x$ the \emph{largest} of
the labels among boxes in $\nu/\lambda$ covered by~$x$. The {\bf reverse
rectification} ${\tt revrectification}(T)$ is the end result of
the iterated application of these slides, ending with a standard
filling of a rotated shape.

\begin{Proposition}
\label{prop:reverse}
Fix a skew shape $\nu/\lambda\subseteq \Lambda_{G/P}$ and 
$T\in {\rm SYT}_{G/P}(\nu/\lambda)$. 
\begin{itemize}
\item[(a)] The analogue for reverse jeu de taquin of Theorem~\ref{thm:proctor}
holds.
\item[(b)] Suppose $y\in \Lambda_{G/P}$ is vacated by
${\tt jdt}_{x}(T)$, then ${\tt revjdt}_{y}({\tt jdt}_{x}(T))=T$.
\item[(c)] Suppose $y\in \Lambda_{G/P}$ is vacated by
${\tt revjdt}_{x}(T)$, then ${\tt jdt}_{y}({\tt revjdt}_{x}(T))=T$.
\end{itemize} 
\end{Proposition}
\begin{proof}
Since $\Lambda_{G/P}$ is self-dual, reverse jeu de taquin slides
is identified with ordinary jeu de taquin slides where the
labels ``$i$'' play the usual role of the labels ``$|\nu/\lambda|-i+1$''.
This reduces (a) to Theorem~\ref{thm:proctor}. The claims (b) and (c)
are easy inductions on the size of $T$.
\end{proof}

\subsection{The infusion involution}
Given $U\in {\rm SYT}_{G/P}(\nu/\lambda)$, a specific choice of jeu de
taquin slides usable to rectify $U$ can be recorded as a tableau $T\in 
{\rm SYT}_{G/P}(\lambda)$. 
Suppose $U$ rectifies to $X$ with ${\tt shape}(X)=\gamma$.  It turns
out that there
is a natural choice of 
tableau $Y$ of shape $\nu/\gamma$ which rectifies to $T$.  In fact, the
map taking $(U,T)$ to $(X,Y)$, which we will call {\tt infusion}, and which
we define below, is an involution.  

Given $T\in {\rm SYT}_{G/P}(\lambda)$ and $U\in {\rm SYT}_{G/P}(\nu/\lambda)$,
we define ${\tt infusion}(T,U)$ to be a pair of tableau $(X,Y)$ with
$X\in {\rm SYT}_{G/P}(\gamma)$ and $Y\in {\rm SYT}_{G/P}(\nu/\gamma)$
(for some $\gamma\in {\mathbb Y}_{G/P}$ with $|\gamma|=|\nu/\lambda|$) 
as follows: place $T$ and $U$
inside $\Lambda_{G/P}$ according to their shapes. 
Now remove the largest label ``$m$'' that appears in $T$,
say at box $x\in\lambda$. Since $x$ necessarily lies next to 
$\nu/\lambda$, apply the slide ${\tt jdt}_{x}(U)$,
leaving a ``hole'' at the other side of $\nu/\lambda$.
Place ``$m$'' in that hole and repeat moving the labels
originally from $U$ until all labels of $T$
are exhausted. In particular, we declare that the labels placed in the created holes
at each step never move for the duration of the procedure. The resulting
straight shape tableau of shape $\gamma$ and skew
tableau of shape $\nu/\gamma$ are $X$ and $Y$ respectively.

\begin{Example}
Let $G/P=Gr(3,{\mathbb C}^7)$, $\lambda=(2,1)$ and $\nu=(3,3,2)$.
The elements of $T$ and $U$ are depicted below, with the labels of
the former are underlined and the labels of the latter 
are given in bold. We also compute ${\tt infusion}(T,U)$ as a sequence
of jeu de taquin slides, where at each
stage the labels of (the eventual) $Y$ are marked with ``$\star$'':
\[(T,U)=\tableau{{{\bf 4} }&{{\bf 5}}&{\ }&{\ }\\
{\underline 3 }&{{\bf 2}}&{{\bf 3} }&{\ }\\
{\underline 1 }&{{\underline 2}}&{{\bf 1} }&{\ }
}\mapsto
\tableau{{{\bf 4} }&{{\bf 5}}&{\ }&{\ }\\
{\bf 2 }&{\bf 3}&{3\star}&{\ }\\
{\underline 1 }&{\underline 2}&{{\bf 1} }&{\ }
}\mapsto
\tableau{{{\bf 4} }&{{\bf 5}}&{\ }&{\ }\\
{\bf 2 }&{\bf 3}&{3\star}&{\ }\\
{\underline 1 }&{\bf 1}&{2\star }&{\ }
}\mapsto
\tableau{{{\bf 4} }&{1\star}&{\ }&{\ }\\
{\bf 2 }&{\bf 5}&{3\star}&{\ }\\
{\bf 1 }&{\bf 3}&{2\star }&{\ }
}={\tt infusion}(T,U)=(X,Y).
\] 
Hence $\gamma=(3,2)$.
\end{Example}

\begin{Theorem}
\label{thm:infusion}
The procedure {\tt infusion} defines an involution on
\[{\tt LockingSYT}(G/P):=\bigcup_{\nu,\gamma\in {\mathbb Y}_{X}} {\rm SYT}_{G/P}(\gamma)\times 
{\rm SYT}_{G/P}(\nu/\gamma).\]
In particular, let $\lambda,\mu,\nu\in {\mathbb Y}_{G/P}$ and set
\begin{multline}\nonumber
{\mathcal I}_{\lambda,\mu}^{\nu}(G/P)=\{(T,U): {\tt shape}(T)=\lambda,\ 
{\tt shape}(U)=\nu/\lambda, \ {\tt infusion}(T,U)=(X,Y) \\ \mbox{ where
${\tt shape}(X)=\mu$}\}\subseteq {\tt LockingSYT}(G/P).
\end{multline}
Then ${\tt infusion}$ bijects ${\mathcal I}_{\lambda,\mu}^{\nu}(G/P)$ and
${\mathcal I}_{\mu,\lambda}^{\nu}(G/P)$.
\end{Theorem}
\begin{proof}
Consider the procedure {\tt revinfusion}: given 
$T\in {\rm SYT}_{G/P}(\lambda)$, $U\in {\rm SYT}_{G/P}(\nu/\lambda)$ placed
inside $\Lambda_{G/P}$ as in the above description of {\tt infusion},
instead remove the smallest label ``$s$'' that appears in $U$,
say at box $x\in \nu/\lambda$. Apply ${\tt revjdt}_{x}(T)$,
leaving a ``hole'' at the other side of $\lambda$. Place ``$s$''
in that hole and repeat until all labels of $U$ are used.
This is indeed the inverse map of {\tt infusion}, as seen by inductively
applying Proposition~\ref{prop:reverse}(b), 
i.e., ${\tt revinfusion}({\tt infusion}(T,U))=(T,U)$.

Thus, the
``involution'' assertion of the theorem amounts to showing that \linebreak
${\tt revinfusion}(T,U)={\tt infusion}(T,U)$ for all $(T,U)$. 
View both procedures and the jeu de taquin
slides as a sequence of ``swaps'' of
the labels of adjacent boxes in $\Lambda_{G/P}$. We prove 
that for each label in $T$ and $U$, 
those swaps involving the label are the same (i.e., 
in the same order and in the same
position inside $\Lambda_{G/P}$,
although possibly occurring at different times in the
overall swap sequence). Then since the path taken by any label
is the same in both procedures, the results are the same. In fact,
it suffices to establish that $T$'s labels undergo the same
swaps in both operations, as the claim about $U$'s labels is then 
implied.

It is easy to check from the definitions of {\tt infusion} and 
{\tt revinfusion} that the swaps
involving the largest label ``$m$'' of $T$ are the same in both procedures.
Since after {\tt infusion} moves ``$m$''
into its final position in $\Lambda_{G/P}$ it never moves again, for
$|\lambda|\geq 2$ the
remainder of the {\tt infusion} procedure is an application of
{\tt infusion} to $({\widetilde T}, {\widetilde U})$ where ${\widetilde T}$
is $T$ with ``$m$'' removed and ${\widetilde U}$ is $U$ after ``$m$''
has moved through it. By induction, ${\tt revinfusion}({\widetilde T}, {\widetilde U})={\tt infusion}({\widetilde T}, {\widetilde U})$ and all swaps 
involving labels of ${\widetilde T}$ and ${\widetilde U}$ are the same.

It remains to show that the swaps not using ``$m$'' are the same
in both ${\tt revinfusion}(T,U)$ and
${\tt revinfusion}({\widetilde T}, {\widetilde U})$. Define two
sequences $T_0:=T, T_1,\ldots, T_{|\nu/\lambda|}$ and ${\widetilde T}_0:={\widetilde T}, 
{\widetilde T}_1,
\ldots {\widetilde T}_{|\nu/\lambda|}$,
where $T_i$ and ${\widetilde T}_{i}$ are the tableaux resulting from moving the labels ``$1$''
through ``$i$'' of $U$ (respectively ${\widetilde U}$) through $T$ during 
${\tt revinfusion}(T,U)$ (respectively 
${\tt revinfusion}({\widetilde T}, {\widetilde U})$). Similarly define the
pair of sequences
$U_0:=U, U_1,\ldots, U_{|\nu/\lambda|}$ and ${\widetilde U}_0:={\widetilde U}, {\widetilde U}_1,
\ldots {\widetilde U}_{|\nu/\lambda|}$,
which are derived from $U$ and ${\widetilde U}$ respectively after the said moving
of labels.

We show by induction on $i\geq 0$ that $T_i$ is ${\widetilde T}_{i}$
with an added corner box containing ``$m$'' and ${\widetilde U}_i$ is obtained by applying to $U_i$
a jeu de taquin slide into that box occupied by ``$m$''. 

The base case $i=0$
holds by construction. For the induction step, there are two cases to consider.
If the ``$m$'' in $T_{i-1}$ is not adjacent to the ``$i$'' in $U_{i-1}$, then
``$i$'' occupies the same (corner) box in both $U_{i-1}$ and
${\widetilde U}_{i-1}$ since the jeu de taquin slide that makes up the difference
between these two tableau does not affect the position of the label ``$i$''.
Therefore the same moves will be made as we pass ``$i$''
through $T_{i-1}$ and ${\widetilde T_{i-1}}$ (in particular, the label ``$m$'' never moves). 
Hence the desired conclusion for $T_i$ and ${\widetilde T}_i$ holds.
In addition, notice $U_i$ and ${\widetilde U}_i$ only differ from 
$U_{i-1}$ and ${\widetilde U}_{i-1}$ respectively by removing the label ``$i$''.
Moreover, the jeu de taquin slide of $U_{i-1}$ into the box labeled ``$m$'' has the same
effect 
as the jeu de taquin slide of $U_i$ into that box. Thus, $U_i$ and 
${\widetilde U}_i$ differ by a jeu de taquin slide into that box, as desired.

Otherwise, if the ``$m$'' in $T_{i-1}$ is adjacent to the ``$i$'' in $U_{i-1}$, then to
obtain $T_i$, the first swap that occurs is between the ``$m$'' and the ``$i$''. But after
this initial swap, the ``$i$'' will be in the same location as the ``$i$'' in ${\widetilde U}_{i-1}$
and so the same swaps will be made as we pass to $T_i$ and ${\widetilde T}_{i}$, so 
these latter two tableau satisfy the conclusion. Also, since the jeu de taquin slide
that makes the difference between $U_{i-1}$ and ${\widetilde U}_{i-1}$ involves 
swapping ``$i$'' and ``$m$'' as the first step, it is clear 
that $U_i$ and ${\widetilde U}_i$ satisfy the conclusion. This
completes the induction argument, and hence the theorem.
\end{proof}

\begin{Corollary}
\label{cor:proctor.implicit}
If $\mu\!\in\! {\mathbb Y}_{G/P}$ and $U\!\in\! {\rm SYT}_{G/P}(\mu)$ then
$\#\{W\!\!\in\! {\rm SYT}(\nu/\lambda)\! :\! {\tt rectification}(W)\!=U\}$
is independent of the choice of $U\in {\rm SYT}_{G/P}(\mu)$.
\end{Corollary}
\begin{proof}
First consider ${\mathcal I}_{\mu,\lambda}^{\nu}(G/P)$.
Since choosing $U\in {\rm SYT}_{G/P}(\mu)$ amounts to 
a choice of sequence of jeu de taquin slides that
rectify $V\in {\rm SYT}_{G/P}(\nu/\mu)$, and by Theorem~\ref{thm:proctor}, 
any choice leads to the same rectification, we conclude that the number of
times a particular tableau $U\in {\rm SYT}_{G/P}(\mu)$ appears as the first
component of a pair in ${\mathcal I}_{\mu,\lambda}^{\nu}(G/P)$ is independent
of~$U$. Specifically, it equals  
$\#\{V\in {\rm SYT}_{G/P}(\nu/\mu):
{\tt shape}({\tt rectification}(V))=\lambda\}.$
On the other hand, 
the number of elements of ${\mathcal I}_{\lambda,\mu}^{\nu}(G/P)$ which 
are carried by ${\tt infusion}$ to a pair of tableaux beginning with $U$ 
is 
$f^{\lambda}(G/P)\cdot 
\#\{W\in {\rm SYT}_{G/P}(\nu/\lambda): {\tt rectification}(W)=U\}$.
By Theorem~\ref{thm:infusion}, these two numbers are equal, so:
\begin{multline}\nonumber
\#\{W\in {\rm SYT}_{G/P}(\nu/\lambda): {\tt rectification}(W)=U\}=\\
\frac{\#\{V\in {\rm SYT}_{G/P}(\nu/\mu): {\tt shape}({\tt rectification}(V))=
\lambda\}}{f^{\lambda}(G/P)}
\end{multline}
is independent of~$U$.
\end{proof}

\begin{Proposition}\label{prop:unrot} 
Fix $\mu\in {\mathbb Y}_{G/P}$. The procedures
${\tt rectification}$ and ${\tt revrectification}$ are
mutually inverse bijections between
${\rm SYT}_{G/P}({\tt rotate}(\mu))\ 
\widetilde{\leftrightarrow} \ {\rm SYT}_{G/P}(\mu)$.
\end{Proposition}

\begin{proof} 
By Theorem~\ref{thm:infusion} specialized to $\nu=\Lambda_{G/P}$, {\tt
rectification} and {\tt revrectification} are mutually inverse bijections
between the set of all fillings of rotated shapes and the set of all
fillings of straight shapes.  Given $T\in {\rm SYT}_{G/P}(\mu)$, let ${\tt
shape}({\tt revrectification}(T))={\tt rotate} (\alpha)$ for some
$\alpha\in {\mathbb Y}_{G/P}$.  By Corollary~\ref{cor:proctor.implicit},
every tableau in ${\rm SYT}_{G/P}(\mu)$ appears as the rectification of
some filling of shape ${\tt rotate(\alpha)}$.  By
Corollary~\ref{cor:proctor.implicit} applied to \linebreak {\tt revrectification},
any filling of shape ${\tt rotate}(\alpha)$ can be obtained by
reverse-rectifying some filling of shape $\mu$, so any filling of shape
${\tt rotate}(\alpha)$ rectifies to a filling of shape $\mu$.  Thus, {\tt
rectification} and {\tt revrectification} are mutually inverse bijections
between fillings of shape ${\tt rotate}(\alpha)$ and of shape $\mu$.  We
can therefore define a bijection $\Psi$ which takes a shape $\alpha$ to
the unique shape of the rectification of any standard filling of ${\tt
rotate}(\alpha)$.  We now show $\Psi$ is the identity map.

First, we show that $\Psi$ is an automorphism of ${\mathbb Y}_{G/P}$.  
Suppose $\alpha$ consists of a shape $\beta$ together with a single
additional 
corner box $x$.  It suffices to show that
$\Psi(\beta)\subseteq\Psi(\alpha)$.  
Fix a filling $T$ of ${\tt rotate}(\alpha)$
in which ${\tt rotate}(x)$ has label ``$1$''.  Let 
the tableau
${\widetilde T}$ be $T$ with the box labeled 1 removed. 
Clearly jeu de taquin methods apply to ${\widetilde T}$.
Pick any sequence of jeu de taquin slides
$T_0 = T, \ T_1={\tt jdt}_{x_1}(T_0),\ \ldots,\  T_i = {\tt jdt}_{x_i}(T_{i-1}),\ \ldots,\  {\tt rectification}(T)$
rectifying $T$. We can define a parallel sequence of jeu de taquin slides
$\{{\widetilde T}_i={\tt jdt}_{{\widetilde x}_{i}}({\widetilde T}_{i-1})\}$
starting with ${\widetilde T}_0 ={\widetilde T}$ where ${\widetilde x}_i
= x_i$ if the label of $T_{i-1}$ moving into $x_i$ is not ``$1$'' and
${\widetilde x}_{i}$ is the box in $T_{i-1}$ with label ``$1$'' otherwise.
It is easy to check that for each $i$, $T_i$ with the ``$1$'' removed
is ${\widetilde T}_{i}$. Also, when $T=T_k$ has been rectified,
${\widetilde T}_k$ is a filling of $\Psi(\alpha)/(1)$.
Since a jeu de taquin slide removes an
``outside corner'' from a skew shape, then $\Psi(\alpha)$ 
is 
$\Psi(\beta)$ with such a corner added.  
Hence, $\Psi$ takes covering relations to covering
relations, as desired.

Let $Q$ be an arbitrary finite poset, and let $D(Q)$ be the lattice of
order ideals (down-closed sets) in $Q$.  Let $\Aut(Q)$ and 
$\Aut(D(Q))$ denote the group of poset automorphisms of $Q$ and $D(Q)$, respectively.
There is a natural inclusion from $\Aut(Q)$ to $\Aut(D(Q))$.  This 
inclusion is actually a group isomorphism: from $\phi
\in\Aut(D(Q))$, the corresponding automorphism of $Q$ can be recovered
by restricting $\phi$ to the principal order ideals of 
$D(Q)$ (which form a poset canonically isomorphic to $Q$). 

So $\Psi$ induces a poset automorphism of
$\Lambda_{G/P}$.  For $B_n$ ($n\geq 4$),
$C_n$, $D_n$ (for the Orthogonal Grassmannians), $E_6$ and $E_7$
the only poset automorphism is the identity. Hence $\Psi$ is 
the identity in these cases. For the remaining
cases, check that $\Psi$ is
the identity on principal lower order ideals of $\Lambda_{G/P}$. This
is straightforward.
\end{proof}

Let $\{e_{\lambda,\mu}^{\nu}(X)\}$ be computed
by the rule of the Main Theorem.

\begin{Corollary}
\label{cor:symmetry}
$e_{\lambda,\mu}^{\nu}(X)=e_{\mu,\lambda}^\nu(X)
=e_{{\tt rotate}(\nu^c),\mu}^{{\tt rotate}(\lambda^c)}(X)$.
\end{Corollary}
\begin{proof}
By Theorem~\ref{thm:infusion} and Corollary~\ref{cor:proctor.implicit}, we have
$e_{\lambda,\mu}^{\nu}(X)=|{\mathcal I}_{\lambda,\mu}^{\nu}(X)|/f^\lambda(X) 
f^\mu(X)=$ \linebreak $|{\mathcal I}_{\mu,\lambda}^{\nu}(X)|/f^\mu(X) f^\lambda(X)=e_{\mu,\lambda}^{\nu}(X).$
For the remaining equality, fix $U\in {\rm SYT}_{X}({\tt rotate}(\mu))$ 
and note that since ${\tt rotate}(\nu/\lambda)={\tt rotate}(\lambda^c)/{\tt rotate}(\nu^c)$, by
Proposition~\ref{prop:reverse}:
\begin{equation}
\label{eqn:symmetryset}
e_{{\tt rotate}(\nu^c),\mu}^{{\tt rotate}(\lambda^c)}(X)=\#\{T\in {\rm SYT}_{X}(\nu/\lambda): {\tt revrectification}(T)=U\}.
\end{equation}
By Theorem~\ref{thm:proctor} and Proposition~\ref{prop:reverse}(b)
combined, it follows that each $T\in {\rm SYT}_{X}(\nu/\lambda)$ 
in~(\ref{eqn:symmetryset}) satisfies
 ${\tt rectification}(T)={\tt rectification}(U)$.
By Proposition~\ref{prop:unrot}, \linebreak
${\tt shape}({\tt rectification}(U))=\mu$.
Hence $e_{\lambda,\mu}^{\nu}(X)\geq 
e_{{\tt rotate}(\nu^c),\mu}^{{\tt rotate}(\lambda^c)}(X)$. 
Reversing the roles of $\lambda$ and $\nu$ above gives equality. 
\end{proof}

\section{Proof of the main theorem}

\subsection{Schubert like numbers}
Fix a cominuscule flag variety $X=G/P$ and 
a collection of real numbers $\{d_{\lambda,\mu}^{\nu}(X)\}$.
It is useful to call $\{d_{\lambda,\mu}^{\nu}(X)\}$ {\bf Schubert like} if (I)-(IV)
below hold:
\begin{itemize}
\item[(I)] (\emph{$S_3$-symmetry}) 
$d_{\lambda,\mu}^{\nu}(X)=d_{\mu,\lambda}^{\nu}(X)=
d_{{\tt rotate}(\nu^c),\mu}^{{\tt rotate}(\lambda^c)}(X)=d_{\lambda,{\tt rotate}(\nu^c)}^{{\tt rotate}(\mu^c)}(X)=d_{{\tt rotate}(\nu^c),\lambda}^{{\tt rotate}(\mu^c)}(X)$\linebreak
$=d_{\mu,{\tt rotate}(\nu^c)}^{{\tt rotate}(\lambda^c)}(X)$; 
\item[(II)] (\emph{Codimension}) $d_{\lambda,\mu}^{\nu}(X)=0$ unless
$|\lambda|+|\mu|=|\nu|$;
\item[(III)] (\emph{Containment}) 
If $\lambda\not\subseteq\nu$ then 
$d_{\lambda,\mu}^{\nu}(X)=0$ for all $\mu\in {\mathbb Y}_{X}$; and 
\item[(IV)] (\emph{Iterated box product}) 
$\sum_{|\gamma|=|\nu/\lambda|} f^\gamma(X) d_{\lambda,\gamma}^\nu(X) 
2^{{\tt shortroots}(\gamma)-{\tt shortroots}(\nu/\lambda)}=f^{\nu/\lambda}(X)$,
where
$f^{\nu/\lambda}(X)=|{\rm SYT}_{X}(\nu/\lambda)|$ and 
$f^{\gamma}(X)=|{\rm SYT}_{X}(\gamma)|$. 
\end{itemize}

\begin{Proposition}
\label{prop:SchubnumsSchublike}
For any cominuscule flag variety $X$,
$\{c_{\lambda,\mu}^{\nu}(X)\}$ is Schubert like.
\end{Proposition}

The Monk-Chevalley formula states that for any $\beta\in\Delta$ and $w\in W$
\begin{equation}
\label{eqn:monk}
\sigma_{s_{\beta}}\cdot\sigma_{w}=\sum_{\alpha\in\Phi^{+}, \ \ell(ws_{\alpha})=\ell(w)+1} n_{\alpha\beta}\frac{(\beta,\beta)}{(\alpha,\alpha)} \ \sigma_{ws_{\alpha}}.
\end{equation}
Here $n_{\alpha\beta}$ is the coefficient of $\beta$ in the expansion of 
$\alpha$ into simple roots, and $(\bullet,\bullet)$ is the inner product
defined on the span of $\Delta$, as determined by the Cartan matrix.

\begin{Lemma}
\label{lemma:grassobs}
Let $\beta=\beta(P)$ and $w\in W$ be Grassmannian at $\beta(P)$. Then in
(\ref{eqn:monk})
\begin{itemize}
\item[(i)] $n_{\alpha,\beta}\in\{0,1\}$
\item[(ii)] $ws_{\alpha}$ is Grassmannian at position $\beta(P)$ whenever
$n_{\alpha\beta}\neq 0$
\end{itemize}
\end{Lemma}
\begin{proof}
(i) holds by the definition of $\beta=\beta(P)$ being cominuscule.
(ii) uses our discussion about 
$H^{\star}(G/P)\hookrightarrow H^{\star}(G/B)$
in the paragraph before Proposition~\ref{prop:lower}. 
\end{proof}

\begin{Proposition}
\label{prop:box_mult}
For any cominuscule flag variety $X$,
\begin{equation}
\label{eqn:comin_monk}
\sigma_{\Box}\cdot\sigma_{\lambda}=\sum_{\mu\in{\mathbb Y}_{X} 
\mbox{ \small {\rm and }}\mu/\lambda \mbox{ \small 
 {\rm is a single box}}} 2^{{\tt shortroots}(\mu/\lambda)}\sigma_{\mu}
\end{equation}
Moreover, $\sigma_{\Box}^{i}=\sum_{|\gamma|=i} 
f^{\gamma}(X) 2^{{\tt shortroots}(\gamma)} \sigma_{\gamma}$.
\end{Proposition}
\begin{proof}
The first claim is proved by comparing (\ref{eqn:monk}) and
(\ref{eqn:comin_monk}) and applying Lemmas~\ref{lemma:bruhatweak} and~\ref{lemma:grassobs}.
The second holds by the definition of $f^{\gamma}(X)$, 
since each standard tableau is inductively built up by adding the box 
labeled ``$k$'' at the $k^{th}$ step.
\end{proof}

\noindent
\emph{Proof of Proposition~\ref{prop:SchubnumsSchublike}:}
(I) and (II) hold by the geometric definition of $c_{\lambda,\mu}^{\nu}(X)$
while (III) holds by Corollary~\ref{cor:overlap}. For (IV), using 
the definition of $f^{\nu/\lambda}(X)$, 
Proposition~\ref{prop:box_mult}, and Corollary~\ref{cor:overlap} we have 
\[\sigma_{\Lambda_{X}}f^{\nu/\lambda}(X)2^{{\tt shortroots}(\nu/\lambda)}=\sigma_{\lambda}\sigma_{\Box}^{|\nu/\lambda|}\sigma_{{\tt rotate}(\nu^c)}=\sum_{|\gamma|=|\nu/\lambda|} 
f^{\gamma}(X) c_{\lambda,\gamma}^{\nu}(X) 2^{{\tt shortroots}(\gamma)}\sigma_{\Lambda_{X}}.\qed \]
\subsection{Cominuscule recursions}
Let ${\widetilde X}={\widetilde G}/{\widetilde P}$ be a second cominuscule
flag variety. Define a {\bf cominuscule recursion} to be a poset injection 
$\Theta:\Lambda_{\widetilde X}\hookrightarrow \Lambda_{X}$
such that $\Lambda_{X}$ is a disjoint union of $\Theta(\Lambda_{\widetilde X})$
 and of $L(\Theta)$ and 
$\Gamma(\Theta)$, which are 
subsets of $\Lambda_X$ whose elements are all either incomparable with or
below (respectively, above) every element of $\Theta(\Lambda_{\widetilde X})$.
Clearly:
\begin{Definition-Lemma}
If $\lambda\in {\mathbb Y}_{X}$ then 
${\overline \lambda}:=\Theta^{-1}(\lambda)$ is in 
${\mathbb Y}_{\widetilde X}$. Also, if  
$\gamma\in {\mathbb Y}_{\widetilde X}$ then 
${\hat \gamma}:=\Theta(\gamma)\cup L(\Theta)$ is in ${\mathbb Y}_X$.
\end{Definition-Lemma}

Fix $\lambda,\mu,\nu\in {\mathbb Y}_{X}$ such that
\begin{equation}
\label{eqn:basic_redux.assumption}
\lambda\subseteq\nu, L(\Theta)\subseteq \lambda \mbox{\  and \ }
\Gamma(\Theta)\subseteq \nu^c.
\end{equation} 
Then $d_{\lambda,\mu}^{\nu}(X)$ is {\bf ${\Theta}$-recursive} if
\begin{equation}
\label{eqn:basic_redux}
d_{\lambda,\mu}^{\nu}(X)=\sum_{\gamma\in {\mathbb Y}_{\widetilde X}}c_{{\overline\lambda},{\gamma}}^ {\overline\nu}({\widetilde X})\ 
d_{L(\Theta), \mu}^{\hat\gamma}(X).
\end{equation}
A collection $\{d_{\lambda,\mu}^{\nu}(X)\}$ is $\Theta$-recursive
if each $d_{\lambda,\mu}^{\nu}(X)$ is, 
whenever (\ref{eqn:basic_redux.assumption}) holds.

Recall that $\{e_{\lambda,\mu}^{\nu}(X)\}$ are the numbers computed
by the rule of the Main Theorem, cf., Corollary~\ref{cor:symmetry}.

\begin{Theorem}
\label{thm:combo_redux}
Fix a cominuscule recursion $\Theta:\Lambda_{\widetilde X}\to \Lambda_{X}$
and assume $e_{\lambda,\mu}^{\nu}({\widetilde X})=c_{\lambda,\mu}^{\nu}({\widetilde X})$ for all $\lambda,\mu,\nu\in {\mathbb Y}_{\widetilde X}$. Then
$\{e_{\lambda,\mu}^{\nu}(X)\}$ is $\Theta$-recursive. 
\end{Theorem}
\begin{proof}
Construct standard fillings of $\nu/\lambda$ rectifying to a fixed
$S\in {\rm SYT}_{X}(\mu)$ in two steps. First choose 
$\gamma\in{\mathbb Y}_{\widetilde X}$ and one of 
the $e_{L(\Theta),\mu}^{\hat\gamma}(X)$ tableaux
$T\in {\rm SYT}_{X}({\hat\gamma}/L(\Theta))$ 
rectifying to $S$.
By Corollary~\ref{cor:proctor.implicit}, there are 
$c_{{\overline\lambda},\gamma}^{\overline{\nu}}({\widetilde X})$ ways to fill $\nu/\lambda$ 
that rectify to $T$ (i.e., rectify to $T$ viewed as a standard tableau of
$\gamma\in {\mathbb Y}_{\widetilde X}$). As this filling 
of $\nu/\lambda$ rectifies to $S$, 
``$\geq$'' for (\ref{eqn:basic_redux}) holds.

For ``$\leq$'', given a standard filling of $\nu/\lambda$ 
that rectifies to $S$,
we want to show that it can be seen to arise as one the above 
standard fillings.
By Theorem~\ref{thm:proctor}, we can start rectifying 
by exclusively choosing boxes 
$x\in\Theta(\Lambda_{\widetilde X})$ to slide into,
until this is no longer possible. At that point 
we have a tableau in ${\rm SYT}_{X}({\hat\gamma}/L(\Theta))$ for
some $\gamma\in {\mathbb Y}_{\widetilde X}$. 
\end{proof}

We use Theorem~\ref{thm:combo_redux} in the Main Theorem's proof.
Once the latter is proved, we obtain:
\begin{Corollary}
For any cominuscule flag variety $X$ and cominuscule recursion
$\Theta$, $\{c_{\lambda,\mu}^{\nu}(X)\}$ is $\Theta$-recursive.
\end{Corollary}

\subsection{The exceptionals ${\mathbb O}{\mathbb P}^2$ 
and $G_{\omega}({\mathbb O}^3, {\mathbb O}^6)$
\label{sect:E6E7proof}}

\begin{figure}[h]
\begin{center}
\epsfig{file=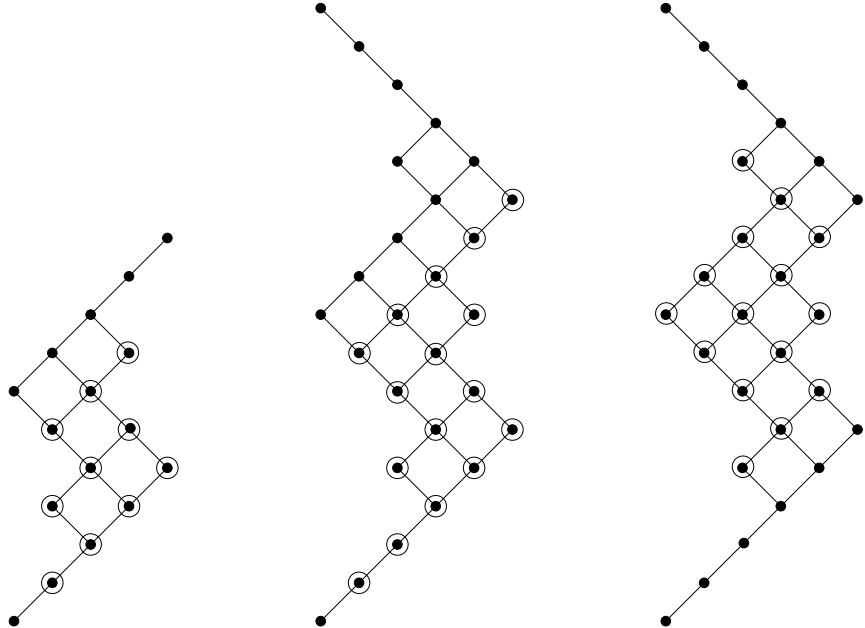, height=8.3cm}
\end{center}
\caption{\label{fig:embeddings} $\Theta_{E_6}, \Theta_{E_7(a)}$ and $\Theta_{E_7(b)}$ respectively (the circled nodes
represent the image of the cominuscule recursion in each case)}
\end{figure}

One has helpful $\Theta$-recursions here:
\begin{itemize}
\item $\Theta_{E_6}: OG(6,12)\to {\mathbb O}{\mathbb P}^2$ 
identifying $\Lambda_{OG(6,12)}$ with
$(1,1,2,3,3,1)/(1)$
\item $\Theta_{E_7(a)}: {\mathbb O}{\mathbb P}^2\to G_{\omega}({\mathbb O}^3, 
{\mathbb O}^6)$ identifying $\Lambda_{{\mathbb O}{\mathbb P}^2}$ with
$(1,1,1,2,4,4,2,1)/(1)$
\item $\Theta_{E_7(b)}: OG(7,14)\to  G_{\omega}({\mathbb O}^3, {\mathbb O}^6)$
identifying $\Lambda_{OG(7,14)}$ with $(1,1,1,2,5,5,3,3)/(1^6)$
\end{itemize} 
Note the ``twist'' in how $\Lambda_{OG(6,12)}$ and $\Lambda_{OG(7,14)}$ are embedded; see Figure~\ref{fig:embeddings}.

The geometric proof of the following proposition is delayed until Section~6.
\begin{Proposition}
\label{prop:geo_arg}
For $X={\mathbb O}{\mathbb P}^2$ and $X=G_{\omega}({\mathbb O}^3, {\mathbb O}^6)$, $\{c_{\lambda,\mu}^{\nu}(X)\}$ is $\Theta_{E_6}$-recursive and respectively
$\Theta_{E_7(a)}$ and $\Theta_{E_7(b)}$-recursive. 
\end{Proposition}

The strategy of the remainder of the proof is that any collection of 
Schubert like numbers $\{d_{\lambda,\mu}^{\nu}(X)\}$ satisfying the 
applicable $\Theta$-recursions above are uniquely determined. Then since
$\{c_{\lambda,\mu}^{\nu}(X)\}$ and $\{e_{\lambda,\mu}^{\nu}(X)\}$ have these
properties, they are the same.

\begin{Lemma}
\label{lemma:uniquedet}
Suppose $\{d_{\lambda,\mu}^{\nu}(X)\}$ is Schubert like. Then
$d_{\lambda,\mu}^{\nu}(X)$ is uniquely determined if any of the following
hold: 
\begin{itemize}
\item[(i)] $X={\mathbb O}{\mathbb P}^2$ and
$d_{\lambda,\mu}^{\nu}(X)$ is $\Theta_{E_6}$-recursive; or 
\item[(ii)] $X=G_{\omega}({\mathbb O}^3, {\mathbb O}^6)$ and
$d_{\lambda,\mu}^{\nu}(X)$ is $\Theta_{E_7(a)}$-recursive; or
\item[(iii)] $X={\mathbb O}{\mathbb P}^2$ or 
$X=G_{\omega}({\mathbb O}^3, {\mathbb O}^6)$ and all but possibly one
$d_{\lambda,\gamma}^{\nu}(X)$ with $|\gamma|=|\nu|-|\lambda|$
is uniquely determined; or
\item[(iv)] $\mu=\emptyset$.
\end{itemize}
\end{Lemma}

\begin{proof}
For (i) and (ii), since $L(\Theta_{E_6})$ and $L(\Theta_{E_7(a)})$ 
is a box,
all $d_{L(\Theta),\mu}^{\hat\gamma}$ 
on the right-hand side of (\ref{eqn:basic_redux})
are determined by (I), (II) and the case $|\nu/\lambda|=1$ of~(IV). 
For (iii), by (IV):
$d_{\lambda,\mu}^{\nu}(X)=(f^{\nu/\lambda}(X)-\sum_{\gamma\neq\mu}f^{\gamma}(X)d_{\lambda,\gamma}^{\nu}(X))/f^{\mu}(X)$
and each of the $d_{\lambda,\gamma}^{\nu}(X)$ is determined, by (i) or (ii). 
Lastly, (iv)
follows from (I) and the case $|\nu/\lambda|=0$ of (IV).
\end{proof}

\begin{Corollary}
\label{cor:uniquedet}
Assume $\{d_{\lambda,\mu}^{\nu}(G_{\omega}({\mathbb O}^3, {\mathbb O}^6))\}$
is Schubert like and $\Theta_{E_7(a)}$ and $\Theta_{E_7(b)}$-recursive. Then 
$d_{\lambda,\mu}^{\nu}(G_{\omega}({\mathbb O}^3, {\mathbb O}^6))$
is uniquely determined if either of the following hold: 
\begin{itemize}
\item[(i)] $|\nu^c|\geq 14$; or 
\item[(ii)] $\lambda=L(\Theta_{E_7(b)})$
and $\Gamma(\Theta_{E_7(b)})\subseteq \nu^c$. 
\end{itemize}
\end{Corollary}
\begin{proof}
For (i), if $|\nu^c|\geq 18$ then 
$\Gamma(\Theta_{E_7(a)})\subseteq\nu^c$, and apply 
Lemma~\ref{lemma:uniquedet}(ii) or (iv). 
If $14\leq |\nu^c|\leq 17$ then all shapes of that size contain
${\tt rotate}(\Gamma(\Theta_{E_7(a)}))$ except one, respectively:
$(1,1,1,2,4,4,1), (1,1,1,2,4,4,2), (1,1,1,2,4,4,2,1)$, $(1,1,1,2,4,4,2,1,1)$, 
and then use \linebreak Lemma~\ref{lemma:uniquedet}(ii) or (iv), or (iii) applied to 
$d_{\lambda,{\tt rotate}(\nu^c)}^{{\tt rotate}(\mu^c)}(X)
=d_{\lambda,\mu}^{\nu}(X)$.

For (ii) if $|\nu^c|\geq 14$ then apply~(i). 
So assume $|\nu^c|\leq 13$. Thus $|\mu|\geq 8$.
If $L(\Theta_{E_7(b)})\subseteq \mu$ then 
$d_{\mu,L(\Theta_{E_7(b)})}^{\nu}$ is $\Theta_{E_7(b)}$-recursive and 
by (\ref{eqn:basic_redux}):
\[
d_{\mu,L(\Theta_{E_7(b)})}^{\nu}(G_{\omega}({\mathbb O}^3, {\mathbb O}^6))
\!=\!\!\!\!\!\!\sum_{\gamma\in {\mathbb Y}_{OG(7,14)}} c_{{\overline\mu},{\gamma}}^ {\overline\nu}(OG(7,14))\ 
d_{L(\Theta_{E_7(b)}),L(\Theta_{E_7(b)})}^{\hat{\gamma}}(G_{\omega}({\mathbb O}^3, {\mathbb O}^6)).
\]
and each nonzero $d_{L(\Theta_{E_7(b)}),L(\Theta_{E_7(b)})}^{\hat{\gamma}}(G_{\omega}({\mathbb O}^3, {\mathbb O}^6))$
is known by (i) since $|{\hat \gamma^c}|=15$.
If \linebreak $L(\Theta_{E_7(b)})\not\subseteq\mu$ then $8\leq |\mu|\leq 10$
and $\mu$ is the only shape of that size
not containing $L(\Theta_{E_7(b)})$. Then $d_{\lambda,\mu}^{\nu}(G_{\omega}({\mathbb O}^3, {\mathbb O}^6))$
is determined by the above argument and Lemma~\ref{lemma:uniquedet}(iii).
\end{proof}

\begin{Proposition}
\label{prop:determined}
Let $X={\mathbb O}{\mathbb P}^2$ or
$X=G_{\omega}({\mathbb O}^3, {\mathbb O}^6)$. If 
$\{d_{\lambda,\mu}^{\nu}(X)\}$ are Schubert like and satisfy the
applicable $\Theta$-recursions, then each $d_{\lambda,\mu}^{\nu}(X)$ is 
uniquely determined 
(and by Propositions~\ref{prop:SchubnumsSchublike} and~\ref{prop:geo_arg} 
thus equal to $c_{\lambda,\mu}^{\nu}(X)$).
\end{Proposition}
\begin{proof}
Suppose $X={\mathbb O}{\mathbb P}^2$. If 
$|\lambda|\!+\! |\mu|\!+\! |\nu^c|< 16
=|\Lambda_{X}|$ then use (II). Otherwise at least one
of $|\lambda|,|\mu|$ or $|\nu^c|$ 
is at least~$6$. By (I) assume it is $|\nu^c|$.
If $|\nu^c|\geq 9$, then 
$\Gamma(\Theta_{E_6})\subseteq \nu^c$ 
and apply Lemma~\ref{lemma:uniquedet}(i) or (iv). 
For $6\leq|\nu^c|\leq 8$ there
is only one shape of that size to which we cannot use 
Lemma~\ref{lemma:uniquedet}(i) or (iv),
namely $(1,1,2,2)$, $(1,1,2,3)$, and $(1,1,2,4)$ respectively. So we can use
Lemma~\ref{lemma:uniquedet}(iii) applied to $d_{\lambda,
{\tt rotate}(\nu^c)}^{{\tt rotate}(\mu^c)}(X)=d_{\lambda,\mu}^{\nu}(X)$.
This completes the proof for this case.

Now let $X=G_{\omega}({\mathbb O}^3, {\mathbb O}^6)$. 
If $|\lambda|+|\mu|+|\nu^c|
<27=|\Lambda_{X}|$ then use (II). Otherwise:
by (I), assume $|\nu^c|\geq {\rm max}(|\lambda|,|\mu|, 9)$, by 
Corollary~\ref{cor:uniquedet}(i) assume
$9\leq |\nu^c|\leq 13$ and by Lemma~\ref{lemma:uniquedet}(ii), assume that 
$d_{\lambda,\mu}^{\nu}(X)$ is
not $\Theta_{E_7(a)}$-recursive. 

Hence, if $|\nu^c|\geq 10$ then
$\Gamma(\Theta_{E_7(b)})\subseteq \nu^c$. Since
$10\leq |\nu^c|\leq 13$, at least one of $\lambda$ or $\mu$ has
size at least 7.
Suppose by (I)
that it is $\lambda$. 
If
$d_{\lambda,\mu}^{\nu}(X)$ is $\Theta_{E_7(b)}$-recursive, 
we use (\ref{eqn:basic_redux}) and 
Corollary~\ref{cor:uniquedet}(ii). 
Note that for the possible values of 
$|\lambda|$,
at most one shape $\gamma$ of that size has
$L(\Theta_{E_7(b)})\not\subseteq \gamma$. 
Thus, if $d_{\lambda,\mu}^{\nu}(X)$ is not
$\Theta_{E_7(b)}$-recursive, it is the unique
$d_{\gamma,\mu}^{\nu}(X)$ with $|\gamma|=|\lambda|$ which is not,
and we use Lemma~\ref{lemma:uniquedet}(iii).

Finally, let $|\nu^c|=9$. Thus $|\lambda|=|\mu|=9$ also.
If $\Gamma(E_7(b))\subseteq\nu^c$ 
then apply the argument of the previous paragraph. Otherwise
${\tt rotate}(\nu^c)=(1,1,1,2,4)$. 
By (I), we are done unless in fact
$\lambda=\mu={\tt rotate}(\nu^c)$. So it remains to consider
$d_{(1,1,1,2,4),(1,1,1,2,4)}^
{{\tt rotate}((1,1,1,2,4)^c)}(X)=d_{(1,1,1,2,4),(1,1,1,2,4)}^{(1,1,1,2,5,4,2,1,1)}$. By (IV): 
$f^{(1,1,1,2,5,4,2,1,1)/(1,1,1,2,4)}(X) = \sum_{|\gamma|=9}f^\gamma(X)
d^{(1,1,1,2,5,4,2,1,1)}_{(1,1,1,2,4),\gamma}(X)$ which in turn equals
$f^{(1,1,1,2,4)}(X) d_{(1,1,1,2,4),(1,1,1,2,4)}^{(1,1,1,2,5,4,2,1,1)}(X) \!+ 
\!\sum_{\gamma\neq (1,1,1,2,4)}\!
f^\gamma(X) d_{(1,1,1,2,4),\gamma}^{(1,1,1,2,5,4,2,1,1)}(X).$
In the latter summation, each $\gamma$ contains 
${\tt rotate}(\Gamma(\Theta_{E_7(b)}))$.
So $d_{(1,1,1,2,4),\gamma}^{(1,1,1,2,5,4,2,1,1)}$ 
is determined by (I) and the argument of the previous paragraph. 
The proposition follows.
\end{proof}

\begin{Proposition}
\label{prop:e_satisfy}
For $X={\mathbb O}{\mathbb P}^2$ and 
$X=G_{\omega}({\mathbb O}^3, {\mathbb O}^6)$,
$\{e_{\lambda,\mu}^{\nu}(X)\}$ are Schubert like, 
satisfy the applicable $\Theta$-recursions, and hence the Main Theorem holds in these 
cases.
\end{Proposition}
\begin{proof}
Clearly (II) and (III) are satisfied. (I) is immediate from 
Corollary~\ref{cor:symmetry}.
For (IV) we need to prove
$\sum_{\gamma,|\gamma|=|\nu/\lambda|}f^{\gamma}(X)e_{\lambda,\gamma}^{\nu}(X)= f^{\nu/\lambda}(X).$
For each of the $f^{\gamma}(X)$ tableaux $T\in {\rm SYT}_{X}(\gamma)$ there are
$e_{\lambda,\gamma}^{\nu}(X)$ tableaux $U\in{\rm SYT}_{X}(\nu/\lambda)$ such
that ${\tt rectification}(U)=T$. This proves ``$\leq$''. 
Conversely, since any $U\in {\rm SYT}_{X}(\nu/\lambda)$
rectifies to some $T\in {\rm SYT}_{X}(\gamma)$, equality holds.
Finally, $e_{\lambda,\mu}^{\nu}(X)$ satisfies the stated recursions
by Theorem~\ref{thm:combo_redux}. The remaining claim follows from
Propositions~\ref{prop:SchubnumsSchublike},~\ref{prop:determined} and the
fact that 
$c_{\lambda,\mu}^{\nu}(OG(6,12))=e_{\lambda,\mu}^{\nu}(OG(6,12))$ and
$c_{\lambda,\mu}^{\nu}(OG(7,14))=e_{\lambda,\mu}^{\nu}(OG(7,14))$, which are
shown in Section~5.5.
\end{proof}

\vspace{-.2in}
\subsection{The quadrics ${\mathbb Q}^{2n-1}$ and ${\mathbb Q}^{2n-2}$}
From Proposition~\ref{prop:box_mult}, it is easy to check
for ${\mathbb Q}^{2n-1}$ (type $B_n$) that
\[
\sigma_{\Box}^{k}=\left\{
\begin{array}{cc}
\sigma_{(1^k)} & \mbox{if $1\leq k<n$,} \\
2\sigma_{(1^k)} & \mbox{otherwise.}
\end{array}
\right\space .
\]
Since $H^{\star}({\mathbb Q}^{2n-1})$ is generated by $\sigma_{\Box}$, the Main Theorem holds in this case.

The case ${\mathbb Q}^{2n-2}$ (type $D_n$), since $\sigma_{\Box}$ does not
generate $H^{\star}({\mathbb Q}^{2n-2})$, we need to also use
Corollary~\ref{cor:overlap}(c) and the following calculations
using Proposition~\ref{prop:box_mult}:
\[
\sigma_{\Box}^{k}=\left\{
\begin{array}{cc}
\sigma_{(1^k)} & \mbox{if $1\leq k\leq n-2$,} \\
\sigma_{(1^{n-3}, 2)}+\sigma_{(1^{n-1})} & \mbox{if $k=n-1$,}\\
2\sigma_{(1^{n-3},2,1)} & \mbox{if $k=n$,}\\
2\sigma_{(1^{n-3},2,2,1^{k-n-1})} & \mbox{otherwise.}\\
\end{array}
\right\space .
\]
and
$\sigma_{\Box}\cdot\sigma_{(1^{n-3},2)}=\sigma_{\Box}\cdot\sigma_{(1^{n-1})}=
\sigma_{(1^{n-3},2,1)}.$

\subsection{Conclusion of the proof of the Main Theorem; the minuscule cases} 
For $Gr(k,{\mathbb C}^n)$, the result
is a mild reformulation of~\cite{Schutzenberger}. 
Similarly, for $LG(n,2n)$ we have restated the work 
of~\cite{Pragacz} and \cite[Theorem~7.2.2]{worley}.
Now, it is known  that the Schubert intersection
numbers for $OG(n+2,2n+4)$ differ from the $LG(n,2n)$ case by a power of $2$
plainly equal to $2^{{\tt shortroots}(\nu/\lambda)
-{\tt shortroots}(\gamma)}$, 
see, e.g.,~\cite[Section~3]{Bergeron.Sottile} and the 
references therein. This, combined with
Proposition~\ref{prop:e_satisfy} proves the ${\mathbb O}{\mathbb P}^2$ and
$G_{\omega}({\mathbb O}^3, {\mathbb O}^6)$ cases. Together with the analysis of
the quadric cases in Section~5.4, this completes the 
cominuscule cases.
For minuscules, the case of $OG(n,2n+1)$ holds by the remarks
in Sections~2.1 and~2.3, while for ${\mathbb P}^{2n-1}$ use 
Proposition~\ref{prop:box_mult}, as done for the odd quadrics. \qed

\section{Cominuscule recursions and Schubert/Richardson variety isomorphisms}

\subsection{Proof of Proposition~\ref{prop:geo_arg}}
Fix lists of cominuscule Lie data 
$(G,B,T(G),P,\Phi(G), W(G))$ and 
$(H, C, T(H), Q, \Phi(H), W(H))$ as in Section~1.1, where 
$T(H)\subseteq T(G)$, $\Phi(H)\subseteq \Phi(G)$,
${\widetilde X}=H/Q$ and $X=G/P$.  Let
${\mathcal Y}_{w}(X), {\mathcal X}_{w}(X)\subseteq X$ respectively be the Schubert
cell and Schubert variety for $wW(G)_P\in W(G)/W(G)_P$.  The 
{\bf opposite Schubert cell} is ${\mathcal Y}^{w}(X):=BwP/P$ and the
{\bf opposite Schubert variety} is 
${\mathcal X}^{w}(X):={\overline{{\mathcal Y}^{w}(X)}}$.
The {\bf Richardson variety} ${\mathcal X}_{u}^{v}(X)$ is the 
reduced and irreducible scheme-theoretic intersection
${\mathcal X}_{u}(X)\cap {\mathcal X}^{v}(X)$.
With obvious adjustments, we use 
this notation for the subvarieties of~${\widetilde X}$.

Below, we only refer to the cominuscule recursions 
$\Theta$ where $\Theta=\Theta_{E_6}, \Theta_{E_7(a)}$ or $\Theta_{E_7(b)}$. 
When unspecified, statements about $\Theta$ refer to all three choices. Also, 
$\beta(P)$ corresponds to the node 1 
of the $E_6$ and $E_7$ Dynkin diagrams from Table 1. Set
$\delta\in W(G)$ to be $s_{\beta_1}$ if 
$\Theta=\Theta_{E_6}, \Theta_{E_7(a)}$
and $s_{\beta_7}s_{\beta_6}s_{\beta_5}s_{\beta_4}s_{\beta_3}s_{\beta_1}$
if $\Theta=\Theta_{E_7(b)}$. 

There is a natural embedding of $W(H)$ into $W(G)$. Discussion of this,
together with the proofs of the following proposition and lemma are
briefly delayed until Section~6.2.

\begin{Proposition}
\label{prop:geomA} 
There is an embedding 
$\eta:{\widetilde X}\hookrightarrow X$ such that:
\begin{itemize}
\item[(I)] $\eta({\mathcal X}^{w}({\widetilde X}))=
{\mathcal X}^{w\delta}_{\delta}(X)$;
\item[(II)] $\eta({\mathcal X}_{w}({\widetilde X}))=
{\mathcal X}^{w_{\rm max}\delta}_{w\delta}(X)$, where
$w_{\max}\in W(H)$ is the maximal length element that is Grassmannian at 
$\beta(Q)$.
\end{itemize}
\end{Proposition}
 
\begin{Lemma}
\label{lem:corresp} Let $w\in W(H)^{Q}$, with
${\mathcal I}(w)=\gamma$.  Then ${\mathcal I}(w\delta)=\hat\gamma$.  
\end{Lemma}

Assuming Proposition~\ref{prop:geomA} and Lemma~\ref{lem:corresp}, we now 
prove Proposition~\ref{prop:geo_arg}.

\begin{Corollary}
\label{cor:geomB}
$\eta({\mathcal X}_{\overline\lambda}^{\overline\nu}({\widetilde X}))
={\mathcal X}_{\lambda}^{\nu}(X)$ and
$\eta_\star([{\mathcal X}_{\overline\lambda}^{\overline\nu}({\widetilde X})])
=[{\mathcal X}_{\lambda}^{\nu}(X)]\in H_{\star}(X,{\mathbb Q})$.
\end{Corollary}
\begin{proof} Let ${\mathcal I}(u)=\overline\lambda$, ${\mathcal I}(v)=\overline\nu$.  
Now,
$\eta({\mathcal X}_{\overline\lambda}^{\overline\nu}
({\widetilde X}))=\eta({\mathcal X}_{u}(\widetilde X)
\cap {\mathcal X}^{v}(\widetilde X))$.
The image is set-theoretically equal to
${\mathcal X}_{u\delta}^{w_{\rm max}\delta}(X)
\cap {\mathcal X}_{\delta}^{v\delta}(X)
=X_{u\delta}(X)\cap X_{\delta}(X)\cap X^{w_{\rm max}\delta}(X)\cap
X^{v\delta}(X)=X^{v\delta}_{u\delta}(X)=X^{\nu}_{\lambda}(X)$. 
Since $\eta$ is an embedding, this is the (scheme-theoretic) image.
The statement about homology follows.
\end{proof}

\begin{Corollary} 
\label{cor:geomC}
$\eta({\mathcal X}^\gamma({\widetilde X}))={\mathcal X}^{{\hat\gamma}}
_{L(\Theta)}(X)$ and 
$\eta_{\star}([{\mathcal X}^{\gamma}({\widetilde X})])=
[{\mathcal X}^{{\hat\gamma}}_{L(\Theta)}(X)]\in H_{\star}(X, {\mathbb Q})$.
\end{Corollary}
\begin{proof} Specialize Corollary~\ref{cor:geomB}: 
$\nu=\hat \gamma$, $\lambda=L(\Theta)$. So  
$\eta({\mathcal X}^{\gamma}(\widetilde X))=
\eta({\mathcal X}^\gamma_\emptyset(\widetilde X))
={\mathcal X}^{\hat\gamma}_{L(\Theta)}(X)$.
\end{proof}

Since Richardson 
varieties are 
homologous to scheme-theoretic unions of Schubert varieties
(or equally, of opposite Schubert varieties), we have:
\begin{equation}
\label{eqn:homologous}
[{\mathcal X}^{\overline \nu}_{\overline \lambda}({\widetilde X})]=
\sum_{\gamma\in{\mathbb Y}_{\widetilde X}} 
c_{{\overline\lambda},\gamma}^{\overline\nu}({\widetilde X})\ [{\mathcal X}^{\gamma}({\widetilde X})]
\in H_{\star}({\widetilde X}, {\mathbb Q}).
\end{equation}
Pushing forward on both sides of (\ref{eqn:homologous})
gives, by Corollaries~\ref{cor:geomB} and~\ref{cor:geomC}:
\begin{equation}
\label{eqn:final_eqn}
[{\mathcal X}^{\nu}_{\lambda}(X)]=
\sum_{\gamma\in{\mathbb Y}_{\widetilde X}} 
c_{{\overline\lambda},\gamma}^{\overline\nu}({\widetilde X}) \ 
[{\mathcal X}^{{\hat\gamma}}_{L(\Theta)}(X)]\in H_{\star}(X,{\mathbb Q}).
\end{equation}
Expanding each $[{\mathcal X}^{{\hat\gamma}}_{L(\Theta)}(X)]$ into 
$\{[X^{\mu}(X)]\}$ and extracting coefficients
on both sides of (\ref{eqn:final_eqn}), we obtain:
$c^\nu_{\lambda,\mu}(X)=\sum_{\gamma\in{\mathbb Y}_{\widetilde X}}
c_{\overline \lambda,\gamma}^{\overline\nu}({\widetilde X})c_{L(\Theta),\mu}^{\hat\gamma}(X)$,
proving Proposition~\ref{prop:geo_arg}.\footnote{It is worthwhile to point
out that these ideas extend to equivariant $K$-theory $K_{T}(X)$.} \qed

\subsection{Proofs of Proposition~\ref{prop:geomA} and 
Lemma~\ref{lem:corresp}\label{subsect:lemma_proofs}}

We fix inclusions of root systems. When $\Theta=\Theta_{E_6}$
the inclusion identifies the nodes $1,2,3,4,5$ of $D_5$
respectively with $6,5,4,3,2$ of $E_6$. Similarly, for 
$\Theta=\Theta_{E_7(a)}$, identify $1,2,3,4,5,6$ of $E_6$ with
$3,2,4,5,6,7$ while for $\Theta=\Theta_{E_7(b)}$ identify
$1,2,3,4,5,6$ of $D_6$ with $1,3,4,5,6,2$ of $E_7$.  
This induces inclusions of the objects of the Lie data for $H$ and $G$;
we assume them for the duration of the paper.
In particular,
$W(H)\subseteq W(G)$ as a 
parabolic subgroup. Let $W(G)_{P}$ be the 
parabolic subgroup of $W(G)$ 
corresponding to omitting node 1, and let $W(H)_{Q}$ be the 
parabolic subgroup of $W(H)\subseteq W(G)$
omitting nodes 1 and 3 when $\Theta=\Theta_{E_6}, \Theta_{E_7(a)}$, 
and omitting nodes 2 and 7 when $\Theta=\Theta_{E_7(b)}$. The following
is a straightforward (finite) check:

\begin{Lemma}\label{lem:psi}
If $\alpha \in \Lambda_{\widetilde X}$ then 
$\Theta(\alpha)=\delta^{-1}\alpha$.  
If $\alpha\in \Phi(H)\setminus (-\Lambda_{\widetilde X})$
then $\delta^{-1}\alpha\not\in \Phi(G)\setminus(-\Lambda_X)$. 
Also ${\mathcal I}(\delta)=L(\Theta)$.\end{Lemma}

\begin{Lemma}\label{lem:refbasic} If $w\in W(H)$, then 
$\ell(w\delta)=\ell(w)+\ell(\delta)$. 
\end{Lemma}
\begin{proof} 
Check that ${\mathcal I}(\delta^{-1})\cap \Phi^+(H)=\emptyset$. 
Since ${\mathcal I}(w)\subseteq \Phi^{+}(H)$, then 
${\mathcal I}(\delta^{-1})\cap {\mathcal I}(w)=\emptyset$. 
Thus if $\alpha\in {\mathcal I}(\delta)$, $\delta(\alpha)\in -{\mathcal I}(\delta^{-1})\subseteq \Phi^-(G)
\setminus \Phi^-(H)\subseteq \Phi^-(G)\setminus -{\mathcal I}(w)$.  Thus, $\alpha
\in {\mathcal I}(w\delta)$ as well. Hence by repeated application 
of~\cite[Lemma~1.6(b)]{Humphreys}, 
$\ell(w\delta)=\ell(w)+\ell(\delta)$.
\end{proof}

\noindent
\emph{ Proof of Lemma~\ref{lem:corresp}.}
Lemma~\ref{lem:refbasic} implies
${\mathcal I}(w\delta)={\mathcal I}(\delta)\cup \delta^{-1} {\mathcal I}(w)$.  
Now apply the first and third parts of Lemma~\ref{lem:psi}.
\qed

\medskip

\begin{Corollary}\label{cor:refC} If $w\in W(H)^{Q}$.
Then $w\delta \in W(G)^{P}$.\end{Corollary}
\begin{proof} By Lemma~\ref{lem:corresp}, ${\mathcal I}(w\delta)\in
{\mathbb Y}_{G/P}$. By Proposition~\ref{prop:lower}, $w\delta\in W(G)^P$.
\end{proof}

\begin{Lemma}\label{lem:ptop} $\delta^{-1}Q\delta\subseteq P$.
\end{Lemma}
\begin{proof} 
Let $U_{\alpha}, \alpha\in\Phi(H)$ denote the \emph{root group}, see, 
e.g.,~\cite[Section~26.3]{Hu2}.
It suffices to show $\delta^{-1}T(H)\delta\subseteq
P$ and $\delta^{-1}U_{\alpha}\delta\subseteq P$
for $\alpha\in \Phi(H)\setminus(-\Lambda_{\widetilde X})$ since 
these subgroups generate $Q$. Now, $\delta^{-1}T(H)\delta\subseteq P$ since
$\delta\in W(G)=N(T(G))/T(G)$.  
By \cite[Theorem 26.3]{Hu2}, $
\delta^{-1}U_\alpha\delta=U_{\delta^{-1}(\alpha)}$. 
By the second part of Lemma~\ref{lem:psi}, 
$\delta^{-1}\alpha\in \Phi(G)\setminus(-\Lambda_{X})$ so
$U_{\delta^{-1}(\alpha)}\subseteq P$.  \end{proof}

\noindent
\emph{Proof of Proposition~\ref{prop:geomA}:}
Pick a representative ${{\check \delta}}\in N(T(G))$ of $\delta$. 
Consider the map $\upsilon: H \to G/P$ defined by the inclusion
of $H$ into $G$, right multiplying by ${\check \delta}$ and 
naturally projecting to $G/P$. This map descends to a well-defined 
set-theoretic map $\eta:H/Q\to G/P$: let $x,y\in H$
with $xQ=yQ$, i.e., $x=yq$ for some $q\in Q$.  Now
$\upsilon(x)=x\delta P/P=yq\delta P/P=y\delta \delta^{-1}q\delta P/P=
y\delta P/P=\upsilon(y)$,
where the second-last equality is by Lemma~\ref{lem:ptop}.

Thus the fibers of $\upsilon$ are unions of cosets $xQ$.
Since $Q$ is a closed subgroup of $H$, the \emph{universal mapping property}
\cite[Theorem 12.1]{Hu2} shows $\upsilon$ 
factors uniquely as a morphism through
$H/Q$ and hence this latter morphism must be equal to $\eta$.

Let $U'_w(H):=\prod_{\alpha\in {\mathcal I}(w^{-1})}U_\alpha(H)\subseteq C$. 
There is a normal form for $H/Q$ is given by:
\begin{equation}
\label{eqn:normalform}
H/Q=\coprod_{w\in W(H)^{Q}} U'_w(H) w Q/Q.\end{equation}
Here we have applied~\cite[Theorem 28.4]{Hu2} together with the
well-known identification of ${\mathcal Y}^w=CwQ/Q\subseteq H/Q$ with 
$CwC/C\subseteq H/C$, 
see, e.g.,~\cite{brion:lectures} (just before Example~1.2.3).

Consider $\eta$ applied to a cell of $H/Q$:
$\eta(U'_w(H) wQ/Q)=U'_w(H) w\delta P/P$. By Corollary~\ref{cor:refC},
$w\delta\in W(G)^P$, so 
it indexes a cell in the decomposition of $G/P$ similar to 
(\ref{eqn:normalform}). The inclusion of $H$ into $G$ embeds 
$U'_w(H)$ into $U'_{w\delta}(G)$   
(as varieties, both being affine spaces) since ${\mathcal I}(w^{-1})\subseteq
{\mathcal I}((w\delta)^{-1})$.
Also,
$H\check\delta\subseteq \overline{C_{-}Q\delta}\subseteq
\overline{B_{-}\delta P}$ (again using Lemma~\ref{lem:ptop}).
Thus
$\eta(CwQ/Q)$ embeds into $Bw\delta P/P \cap \overline
{B_{-}\delta P/P}$.
Taking closures we get
${\mathcal X}^w(H/Q)\subseteq {\mathcal X}^{w\delta}_\delta(G/P)$.
These are irreducible varieties of the same dimension, namely 
$\ell(w)$, so this inclusion is an isomorphism.  

For (II), the argument are similar.  First, $\eta$ embeds
${\mathcal Y}_{w}(\widetilde X)=C_{-}wQ/Q$ 
into $B_{-}w\delta P/P=
{\mathcal Y}_{w\delta}(X).$
Next, note that 
$\eta(H/Q)=\overline{\eta(Cw_{\rm max}Q/Q)}\subseteq 
\overline{Bw_{\rm max}\delta P/P}$.
Thus, $\eta({\mathcal Y}_{w}(\widetilde X))$ embeds into
${\mathcal X}_{w\delta}^{w_{\rm max}\delta}$.
Taking closures, and observing that the varieties on both sides have the
same dimension, namely $\ell(w_{\rm max})-\ell(w)$, we see that the
two varieties are isomorphic.  
\qed

\section{Final remarks and questions; (geometric) representation theory}

\begin{Problem}
Find equivariant, $K$-theoretic and/or quantum analogues of the Main Theorem.
\end{Problem}
This is a standard kind of question in the subject. That being said,
in consultation with Allen Knutson and Terence Tao, we surmise that
there is hope to obtain these generalizations in the cominuscule setting.
For example, in the special case of Grassmannians, \emph{puzzle}
theorems/conjectures generalizing the Littlewood-Richardson
rule exist in each of the three basic directions (as well as some
combinations), see respectively, e.g.,~\cite{knutson.tao:equivariant, 
buch:KLR, buch.kresch.tamvakis} (the latter of which reduces
the quantum problem to the $2$-step flag manifold problem, which is
computed by a conjecture of Knutson). Thus:
\begin{Problem}
\label{problem:second}
Reformulate the Main Theorem via a uniform generalization of the puzzles 
of \cite{knutson.tao:equivariant}.
\end{Problem}
In view of~\cite[Appendix~A]{vakil}, a geometric
motivation for Problem~\ref{problem:second} is that a solution might suggest 
appropriate degenerations of Richardson varieties giving 
a geometric version of the Main Theorem, 
in the spirit of~\cite{vakil}. This geometry could yield
interesting arithmetic consequences, see~\cite{vakil:ind}.

We remark on a representation theoretic interpretation
of the Main Theorem, as told to us by Allen Knutson.
One special property of cominuscule flag manifolds that separates
them from general $G/P$'s is that they are also the only smooth
Schubert varieties on the affine Grassmannian. The geometric Satake 
correspondence of Ginzburg, Mirkovi\'{c}-Vilonen and others relates the
geometry of the affine Grassmannian of $G$ to the representation theory
of the Langlands dual group $G^{\vee}$. This correspondence associates
cominuscule flag manifolds to minuscule representations. 
See, e.g.,~\cite{Mirkovic.Vilonen}. Consequently, the Main Theorem computes
special cases of the natural action of the cohomology of the affine 
Grassmannian on the intersection homology of its Schubert varieties. 
This viewpoint suggests further avenues for potential generalization.

Finally, there is a connection between the cominuscule Schubert intersection numbers
and tensor product multiplicities $m_{\lambda,\mu}^{\nu}(G)$ 
of the finitely many 
irreducible representations of Levi subgroups appearing in the 
exterior algebra of the subspace corresponding to $\Lambda_{G/P}$. 
See~\cite[Section~8]{Kostant} and \cite{Belkale.Kumar}.
\begin{Problem}
When does $c_{\lambda,\mu}^{\nu}(G/P)=m_{\lambda,\mu}^{\nu}(G)$?
\end{Problem}
It is known that the equality holds for the ``single simple factor case'' in
$G=GL_{n}({\mathbb C})$.

\section*{Acknowledgments}

We thank Kevin Purbhoo and Frank Sottile
for inspiring conversations about \cite{Purbhoo.Sottile}. We are most grateful
to Allen Knutson for bringing to our attention 
relations of our work to representation theory and for many other comments,
corrections and insights on this text. 
We would also like to thank Prakash Belkale, Mark Haiman, Joseph Landsberg,
Bernd Sturmfels, Terry Tao and Alexander Woo for
helpful discussions, and John Stembridge for supplying code to compute 
Schur $P-, Q-$ products as well as answering our questions about 
its theory. HT was partially supported by an NSERC Discovery grant. 
AY was partially supported by NSF grant DMS-0601010. This work was partially
completed during a visit by HT to the University of Minnesota, supported
by a McKnight grant awarded to Victor Reiner, and separately,
while AY was an NSERC sponsored visitor at the Fields Institute
in Toronto. Finally, AY would like to thank the Institute for Pure and 
Applied Mathematics at UCLA, where exceptional hospitality
greatly facilitated the writing of this paper 
during an NSF supported visit of April-June, 2006.


\begin{thebibliography}{999999999}
\bibitem[BelKu06]{Belkale.Kumar} P.~Belkale and S.~Kumar,
\emph{Eigenvalue problem and a new product in cohomology of flag 
varieties}, Invent.~Math., to appear, 2006. 
\bibitem[BerSo02]{Bergeron.Sottile} N.~Bergeron and F.~Sottile, 
\emph{A Pieri-type formula for isotropic flag manifolds}, 
Trans. Amer. Math. Soc.,~{\bf 354} (2002), no.~7, 2659--2705.  
\bibitem[BiLa00]{Billey.Lakshmibai} S.~Billey and V.~Lakshmibai, \emph{Singular loci
of Schubert varieties}, Progr.~Math.~{\bf 182}(2000), Birkh\"{a}user, Boston.
\bibitem[Bj83]{Bj} A.~Bj\"orner, \emph{Orderings of Coxeter groups}, in 
\emph{Combinatorics and algebra (Boulder, Colo., 1983)}, 175--195, AMS, Providence 
RI, 1984.
\bibitem[BjEdZi90]{Ed} A.~Bj\"orner, P. Edelman, and G. Ziegler, \emph{Hyperplane
arrangements with a lattice of regions}, Discrete Comput. Geom {\bf 5} (1990),
263--288.
\bibitem[Br05]{brion:lectures} M.~Brion, \emph{Lectures on the geometry of flag varieties}, Topics in cohomological studies of algebraic varieties,  33--85, 
Trends Math., Birkh\"{a}user, Basel, 2005.
\bibitem[Bu02]{buch:KLR} A.~Buch, \emph{A Littlewood-Richardson rule for the $K$-theory of Grassmannians}, 
Acta Math.,~{\bf 189} (2002), 37--78.
\bibitem[BuKrTa03]{buch.kresch.tamvakis} A.~Buch, A.~Kresch and H.~Tamvakis, 
\emph{Gromov-Witten invariants on Grassmannians}, J.~Amer.~Math.~Soc.~{\bf 16} (2003), 901--915.
\bibitem[Co05]{Coskun} I.~Coskun, \emph{A Littlewood-Richardson rule for two-step flag
varieties}, preprint 2005.
\bibitem[Hu90]{Humphreys} J.~E.~Humphreys, \emph{Reflection groups and Coxeter
groups}, Cambridge University Press, 1990.
\bibitem[Hu75]{Hu2} \bysame, \emph{Linear algebraic groups}, Springer Verlag,
New York, 1975.
\bibitem[Kn03]{Knutson:alg} A.~Knutson, \emph{A Schubert calculus recurrence from the noncomplex $W$-action on $G/B$}, \textsf{math.CO/0306304}.
\bibitem[KnTa03]{knutson.tao:equivariant} A.~Knutson and T.~Tao, \emph{Puzzles and (equivariant) cohomology
of Grassmannians}, Duke~Math.~J., {\bf 119} (2003), no. 2, 221--260.
\bibitem[KnYo04]{Knutson.Yong} A.~Knutson and A.~Yong, \emph{A formula for $K$-theory truncation
Schubert calculus}, Internat.~Math.~Res.~Notices, {\bf 70}(2004), 3741--3756.
\bibitem[Kog01]{Kogan} M.~Kogan, \emph{RC-graphs and a generalized Littlewood-Richardson rule},
Internat.~Math.~Res.~Notices, {\bf 15}(2001), 765--782.
\bibitem[Kos61]{Kostant} B.~Kostant, \emph{Lie algebra cohomology and
   the generalized Borel-Weil theorem,} Ann. of Math. {\bf (2) 74}(1961), 
329--387.
\bibitem[MiVi99]{Mirkovic.Vilonen} I.~Mirkovi\'{c} and K.~Vilonen, 
\emph{Perverse Sheaves on affine Grassmannians and Langlands duality},
\textsf{math.AG/9911050}.
\bibitem[Pe06]{Perrin} N.~Perrin, \emph{Small resolutions of minuscule Schubert varieties}, \textsf{arXiv:math.AG/0601117}.
\bibitem[Pra91]{Pragacz} P.~Pragacz, \emph{Algebro-geometric applications of Schur S-
and Q- polynomials}, in Topics in invariant theory, Seminaire d'Algebre
Dubreil-Malliavin 1989--1990 (M.-P. Malliavin ed.), Springer Lecture
Notes in Math.~{\bf 1478}, 130--191, Springer, Berlin, 1991.
\bibitem[Pro04]{Proctor} R.~Proctor, \emph{d-Complete posets generalize Young diagrams for the jeu de taquin property}, preprint, 2004, available at
\textsf{http://www.math.unc.edu/Faculty/rap/}
\bibitem[Pu06]{Purbhoo} K.~Purbhoo, 
\emph{Vanishing and nonvanishing criteria in Schubert calculus}, Internat. 
Math.~Res.~Notices, vol.~2006, Article ID 24590, 38 pages, 2006.
\bibitem[PuSo06]{Purbhoo.Sottile} K.~Purbhoo and F.~Sottile, 
\emph{The recursive nature of cominuscule Schubert calculus}, \textsf{arXiv:math.AG/0607669}.
\bibitem[Sa87]{sagan} B.~E.~Sagan, \emph{Shifted tableau, Schur~$Q-$ functions, and a conjecture
of Stanley}, J.~Comb.~Theory, ser.~A {\bf 45} (1987), 62--103.
\bibitem[Sc77]{Schutzenberger} M.-P.~Sch\"{u}tzenberger, \emph{Combinatoire et repr\'{e}sentation du 
groupe sym\'{e}trique} (Actes Table Ronde CNRS, Univ.~Louis-Pasteur Strasbourg, Strasbourg, 1976),  
pp. 59--113. Lecture Notes in Math., Vol. 579, Springer, Berlin, 1977. 
\bibitem[Va05a]{vakil} R.~Vakil, \emph{A geometric Littlewood-Richardson rule} 
(with an appendix joint with A.~Knutson), Annals of Math., to appear, 2005.
\bibitem[Va05b]{vakil:ind} \bysame, \emph{Schubert induction}, Annals of Math., to appear, 2005.
\bibitem[Wo84]{worley} D.~Worley, \emph{A theory of shifted Young tableau}, Ph.~D.~thesis, M.~I.~T., 1984, available at {\tt http://hdl.handle.net/1721.1/15599} .
\bibitem[Yo06]{Yong} A.~Yong, \emph{Maple 7 code to compute Schubert calculus in G/B}, software available at \\
{\tt http://www.math.umn.edu/$\widetilde{\  }$}{\tt ayong/papers.html}
\end{thebibliography}
\end{document}